\documentclass{article}

\usepackage{amsmath,amssymb,amsthm,bm}
\usepackage{graphicx,subcaption}

\usepackage{amsmath,amsfonts,bm}

\def\1{\bm{1}}
\def\rvx{{\mathbf{x}}}
\def\rvy{{\mathbf{y}}}

\def\vzero{{\bm{0}}}
\def\vone{{\bm{1}}}

\def\vtheta{{\bm{\theta}}}
\def\veta{{\bm{\eta}}}
\def\vphi{{\bm{\phi}}}
\def\vxi{{\bm{\xi}}}
\def\va{{\bm{a}}}

\def\vc{{\bm{c}}}
\def\vd{{\bm{d}}}

\def\vf{{\bm{f}}}

\def\vh{{\bm{h}}}

\def\vl{{\bm{l}}}
\def\vm{{\bm{m}}}

\def\vp{{\bm{p}}}
\def\vq{{\bm{q}}}
\def\vr{{\bm{r}}}

\def\vu{{\bm{u}}}

\def\vw{{\bm{w}}}
\def\vx{{\bm{x}}}
\def\vy{{\bm{y}}}
\def\vz{{\bm{z}}}

\def\mA{{\bm{A}}}
\def\mB{{\bm{B}}}

\def\mH{{\bm{H}}}

\def\mL{{\bm{L}}}

\def\mSigma{{\bm{\Sigma}}}

\DeclareMathAlphabet{\mathsfit}{\encodingdefault}{\sfdefault}{m}{sl}
\SetMathAlphabet{\mathsfit}{bold}{\encodingdefault}{\sfdefault}{bx}{n}

\usepackage[normalem]{ulem}
\usepackage{url}
\usepackage{booktabs}
\usepackage{makecell}
\usepackage{xcolor}
\usepackage{algorithm, algorithmic}
\usepackage{tikz}
\usetikzlibrary{shapes,arrows}

\newcommand{\blue}[1]{#1}

\usepackage{natbib}
 \bibpunct[, ]{(}{)}{,}{a}{}{,}%
\usepackage{soul}
\usepackage{authblk}
\usepackage{hyperref}
\usepackage[hmargin=1in]{geometry}

\newtheorem{proposition}{Proposition}

\title{Decision-calibrated prediction sets for robust power system operations}

\author[1]{Akylas Stratigakos}
\author[2]{Honglin Wen}
\author[3]{Elina Spyrou}
\author[4]{Pierre Pinson}
\affil[1]{UCL Energy Institute, University College London}
\affil[2]{Department of Electrical Engineering, Shanghai Jiao Tong University}
\affil[3]{Department of Electrical and Electronic Engineering, Imperial College London}
\affil[4]{Dyson School of Design Engineering,
Imperial College London}

\begin{document}
\maketitle

\begin{abstract}
Robust optimization offers a tractable approach to balance operating costs and reliability in power systems dominated by weather-dependent renewable uncertainty, but its performance depends critically on the uncertainty set.
Data-driven approaches construct uncertainty sets from historical observations or predictive models and often calibrate them for prescribed coverage,
which suffices for the corresponding operational reliability target under robust feasibility.
However, as operating decisions can also satisfy all constraints for realizations outside the uncertainty set,
coverage-based calibration can lead to unnecessarily large sets and higher operating costs.
We propose \emph{decision-calibrated} conditional multivariate prediction sets that serve as uncertainty sets for robust scheduling.
Our two-step framework separates likelihood-based learning of set location and shape from calibration of set size to meet an operational reliability target.
First, we use norm-based representations that capture contextual information and multivariate dependence while admitting tractable robust formulations.
Second, we establish finite-sample control of the expected downstream constraint-violation loss by calibrating the set-size threshold directly against operational violations.
Since reoptimizing operating decisions as the set-size threshold changes can make the loss non-monotone,
we apply conformal risk control to a monotone upper envelope of the loss on a fixed threshold grid.
We evaluate our framework in 15-minute-ahead energy and reserve scheduling with network constraints, which we formulate as a robust DC optimal power flow problem with affine recourse.
Numerical experiments on a modified 5-bus system and the IEEE RTS-GMLC system
show that decision-based calibration tracks prescribed reliability targets more closely and reduces operating costs relative to coverage-based calibration.
Comparisons on the 5-bus system show similar benefits relative to end-to-end learning benchmarks.
We further compare contextual transformations and calibration criteria separately and find that the calibration criterion has a larger effect on downstream performance.
\end{abstract}

{\small\noindent\textbf{Keywords:} Robust optimization, multivariate prediction sets, conformal risk control, conformal prediction, decision-based calibration, power systems\par}

\section{Introduction}

\subsection{Context and Motivation}

In short-term grid scheduling applications (from days to minutes ahead), power system operators (SOs) must balance operating costs and reliability \citep{ferc2021}.
Growing uncertainty from weather-dependent renewables, such as wind and solar, makes this task increasingly challenging.
In principle, multi-stage stochastic programming and chance-constrained approaches provide natural ways to incorporate probabilistic forecasts or scenarios into operational decisions \citep{roald2023power}.
In practice, however, SOs often rely on deterministic scheduling formulations combined with reliability checks under contingency scenarios \citep{Hobbs_SF2}, an approach more similar to robust optimization (RO) than to multi-stage stochastic programming.
RO and, by extension, adaptive RO (ARO) formalize this idea by ensuring that the system has enough flexible capacity to manage all forecast-error realizations within a prescribed \emph{uncertainty set}, possibly after affine or other recourse actions \citep{street2025robustness}.

The success of RO hinges on the choice of uncertainty set, which captures possible uncertainty realizations for which operational feasibility is enforced \citep{bertsimas_hertog2020}.
In data-driven RO, a common approach is to construct uncertainty sets directly from historical observations \citep{shang2017data,goerigk2020data}.
Contextual approaches extend this idea by using predictive models to adapt the set to information available at scheduling, such as weather forecasts \citep{chenreddy2022data,bertsimas2025constructing}.
For multivariate uncertainty, set location represents the central prediction, while set shape captures directional dispersion and dependence among random variables.
For instance, axis-aligned boxes can include simultaneous extremes across several variables because they do not represent dependence.
Conditioning on available information also allows the set to adapt to changes in forecast-error dispersion.

The development of contextual approaches connects data-driven RO to probabilistic energy forecasting \citep{hong2020energy}.
Predictive models can represent uncertainty through marginal quantiles, scenarios, or multivariate conditional distributions.
For RO, \emph{prediction sets} provide a set-valued forecast that can be used directly as an uncertainty set.
Quantile-based prediction intervals are standard in energy forecasting.
For use in RO, multivariate prediction sets should capture dependence, adapt to context, and admit tractable robust reformulations.
We construct such sets as sublevel sets of trained score functions \citep{meng2023scores}.

For a fixed location and shape, a size parameter controls the range of realizations protected against.
A common criterion for choosing this parameter is predictive coverage, i.e., the probability that a future realization lies inside the set.
Conformal calibration can guarantee a prescribed coverage level in finite samples under exchangeability \citep{angelopoulos2023conformal}.
Under robust feasibility, achieving this level suffices for meeting the corresponding operational reliability target \citep{hong2021learning}.
However, operating decisions can also satisfy all constraints for realizations outside the set.
Therefore, calibrating set size to a prescribed coverage level can produce larger sets and higher operating costs than needed to meet the reliability target \citep{bertsimas2021probabilistic}.
This motivates calibrating set size directly against downstream operational constraint violations.

\subsection{Aims and Contributions}

We develop a two-step framework for constructing conditional multivariate prediction sets for robust power system operations.
Throughout, \emph{operational reliability} refers to the probability that all considered operational constraints, namely, reserve and transmission constraints under the scheduled recourse policy, are satisfied when renewable production is realized.
We first learn the conditional prediction-set location and shape, then calibrate a scalar size threshold against downstream constraint violations.
We refer to the resulting sets as \emph{decision-calibrated} prediction sets.
Our contributions are threefold.

First, we formulate conditional norm-based prediction sets that support likelihood-based training and tractable RO.
Building on norm-based score models \citep{braun2025minimum}, we derive the training objectives and robust constraint reformulations for the considered score functions, including a sum of the $1$- and $\infty$-norms.
These representations capture multivariate dependence and contextual information while retaining linear or second-order cone formulations for robust scheduling.

Second, under the assumptions of Proposition~\ref{prop:decision_crc}, we establish finite-sample control of the expected downstream constraint-violation loss by calibrating the set-size threshold directly against operational violations.
Although larger uncertainty sets contain all realizations covered by smaller sets,
reoptimizing the schedule and recourse policy can increase violations for individual realizations.
We address this non-monotonicity by constructing a monotone upper envelope over a fixed threshold grid and applying conformal risk control \citep{angelopoulos2024conformal}.

Third, we evaluate the framework in 15-minute-ahead energy and reserve scheduling with network constraints \citep{vancaelenberg2026deployment}.
Across the tested systems, decision-based calibration generally yields constraint satisfaction closer to the prescribed targets and lower operating costs than coverage-based calibration, which selects sets to contain realized renewable production with a prescribed probability.
On the 5-bus system, we also achieve lower average costs and constraint satisfaction closer to the prescribed targets than end-to-end benchmarks that train uncertainty-set parameters using downstream operational losses.
Further comparisons that vary the contextual transformation and calibration criterion separately show that the latter has a larger effect on downstream performance.

\subsection{Related Work}

Multivariate probabilistic forecasting aims to construct calibrated and sharp predictive distributions \citep{gneiting2007probabilistic} that represent spatio-temporal dependence across uncertain variables and forecasting horizons \citep{hong2020energy}.
In power systems, a popular approach combines marginal predictive distributions with a joint dependence model, e.g., a copula, to generate scenarios \citep{pinson2009probabilistic,carmona2024_joint_granular_model,zhang2025weather}.
\citet{pinson2009probabilistic} uses this approach to generate wind power scenarios,
whereas more recent works generate joint scenarios for load, wind, and solar production \citep{carmona2024_joint_granular_model, zhang2025weather}.
Although these scenario generation methods are useful for stochastic optimization,
they do not produce set-based representations needed for RO formulations.

Prediction sets provide a set-valued representation of uncertainty that can be used directly in RO \citep{meng2023scores}.
Data-driven RO methods construct uncertainty sets from observed realizations, focusing on representations that support tractable robust formulations.
\citet{shang2017data} construct polyhedral sets using kernel-based methods and support vector machines,
whereas \citet{goerigk2020data} use unsupervised learning to construct minimum-volume hyperspheres.
Contextual approaches further enable sets to adapt to available information and account for changes in forecast-error dispersion and dependence \citep{chenreddy2022data}.
In power systems, \citet{shen2024set} examines set-valued regression for wind power curves.
\citet{golestaneh2018ellipsoidal} construct ellipsoidal prediction sets by combining point forecasts with conditional covariance estimation,
whereas \citet{golestaneh2018polyhedral} develop polyhedral prediction sets through convex-hull-based post-processing.
\citet{bertsimas2025constructing} construct contextual sets using regression models and mixed-integer optimization.
\citet{andrianesis2024ensembling} combine wind forecasts and then learn uncertainty sets with coverage guarantees,
while \citet{wasilkoff2023day} construct clustered ellipsoidal sets for wind production.

Once the location and shape of a prediction set are learned, its size can be calibrated to achieve a prescribed predictive coverage.
Conformal prediction (CP) provides distribution-free finite-sample coverage guarantees that a future realization lies inside the prediction set with at least the prescribed probability under data exchangeability \citep{angelopoulos2023conformal}.
In split CP, a predictive model is trained first. Nonconformity scores are computed on held-out \textit{calibration} data, and a {score threshold} parameter is selected using an empirical quantile plus a correction term.
Multivariate CP studies score representations that reduce set volume at a prescribed coverage level.
\citet{braun2025minimum} use norm-based scores to learn minimum-volume multivariate prediction sets, including multi-norm formulations,
whereas \citet{tumu2024multi} optimize convex score functions for multimodal prediction sets.
Power system applications include renewable forecasting and forecast aggregation
\citep{yang2025conformal, moradi2026copula},
while \citet{zhou2026hierarchical} use hierarchical CP to quantify uncertainty in distributed energy resource adoption forecasts.

When a conformal prediction set is used in RO, its coverage guarantee provides a lower bound on the probability of joint operational constraint satisfaction under robust feasibility.
\citet{johnstone2021conformal} develop this connection and construct conformal uncertainty sets, focusing on ellipsoids.
More generally, \citet{hong2021learning} formalize a two-step learning-based RO framework in which the uncertainty-set shape is learned from data and the size is calibrated on a holdout set to obtain feasibility guarantees.
Related power system applications use statistically calibrated or CP-based uncertainty sets for robust unit commitment, reserve sizing, data center operation, and resilient planning \citep{xie2025predict,ling2026dynamic,yang2025conformal,chen2025enhancing}.
The conservativeness of coverage-based calibration \citep{bertsimas2021probabilistic} motivates approaches that use downstream performance to guide uncertainty-set construction.

Decision-focused or end-to-end uncertainty-set learning trains set parameters using the performance of downstream decisions.
Namely, the uncertainty set enters a robust optimization problem, the resulting decision determines the training loss, and the set parameters are updated by differentiating through the solution.
The downstream objective can therefore influence both set shape and size.
\citet{wang2023learning} learn contextual uncertainty sets using downstream objectives and a conditional value-at-risk (CVaR) constraint on operational violations, combining implicit differentiation with an augmented Lagrangian formulation.
\citet{chenreddy2024end} combine downstream task loss with a conditional coverage loss.
\citet{yeh2024end} integrate coverage calibration into end-to-end learning and represent convex sets using partially input-convex neural networks \citep{amos2017input}.
In power systems, \citet{gu2024uncertainty} learn uncertainty sets for adaptive robust economic dispatch,
whereas \citet{mieth2023prescribed} learn contextual sets for robust DCOPF with affine recourse using operating costs and a CVaR-based reliability constraint.
Here, CVaR approximates the probabilistic reliability requirement used during training, while the downstream scheduling problem retains its robust constraints.

Motivated by control of downstream decision reliability, \citet{huconformal} construct prediction sets to minimize the expected worst-case loss of the induced robust decisions.
Their reliability requirement bounds the probability that the realized decision loss exceeds this worst-case loss, and their approach combines end-to-end learning and conformal calibration.

Conformal risk control (CRC) extends coverage calibration to controlling the expected losses for a fixed predictive model.
CRC selects a score threshold to control a bounded monotone loss at a prescribed tolerance, with predictive miscoverage as a special case \citep{angelopoulos2024conformal}.
\citet{yeh2025conformal} consider expected and tail risks of monotone losses and integrate CRC into model training.
\citet{zhou2026calibrating} consider the inverse problem, where they evaluate the expected risk associated with candidate parameters, e.g., a threshold controlling the size of an uncertainty set, and quantify miscoverage--regret trade-offs.
For robust scheduling, the violation loss depends on the decision obtained at each threshold.
A realization outside two nested uncertainty sets may satisfy the schedule obtained from the smaller set and violate a constraint under the schedule obtained from the larger set.
Therefore, the resulting violation loss can increase with the threshold, invalidating the CRC monotonicity assumption.
\citet{angelopoulos2026conformalriskcontrolnonmonotonic} establishes risk guarantees for non-monotone losses under stability conditions on the parameter-selection procedure.

Empirical constraint violations can also be used to tune robust decisions without conformal calibration.
\citet{hou2020chance,hou2021data} iteratively adjust a scalar safety parameter using sample-based evaluations of joint constraint violations.
Their OPF application tunes an ellipsoidal uncertainty set under fixed generator participation factors.
The certification framework of \citet{hou2021data} assumes \emph{i.i.d.} uncertainty observations and separates tuning from verification.
An independent verification sample provides a high-confidence feasibility certificate when the verification criterion is satisfied.

We retain the two-step learning and calibration framework of \citet{hong2021learning} and calibrate set size against operational constraint violations.
We first learn the conditional location and shape using a likelihood-based objective, then fix the trained model and calibrate a scalar threshold through robust scheduling.
This confines downstream optimization to calibration and avoids solving and differentiating through robust scheduling problems during predictive model training.
The calibration evaluations can also be reused across reliability targets.
Applying CRC to a monotone upper envelope yields finite-sample control of the expected downstream violation loss under Proposition~\ref{prop:decision_crc}.
The guarantee is marginal over calibration and test observations,
which differs from independent verification methods that provide feasibility certificates when their verification criterion is satisfied \citep{hou2021data}.

\subsection{Paper Structure}

The rest of the paper is organized as follows.
First, we present the problem setup (Section~\ref{sec:setup}) and develop norm-based score models (Section~\ref{sec:learning}).
Next, we formulate the motivating grid scheduling application (Section~\ref{sec:dcopf}), develop the proposed calibration method (Section~\ref{sec:calibration}),
present numerical results (Section~\ref{sec:results}),
and provide conclusions (Section~\ref{sec:conclusions}).

\section{Problem Setup}\label{sec:setup}

In this section, we provide preliminary background (Subsection~\ref{subsec:preliminaries}),
and outline the proposed method (Subsection~\ref{subsec:outline}).

\subsection{Preliminaries}\label{subsec:preliminaries}

Let $\rvy \in \mathcal{Y} \subseteq \mathbb{R}^{d_\rvy}$ be a multivariate random variable, e.g., renewable production across wind farms,
and let $\rvx \in \mathcal{X} \subseteq \mathbb{R}^{d_\rvx}$ denote associated contextual information or \emph{features}, e.g., weather forecasts or historical production.
For normalized renewable production, we take the physical support to be $\mathcal{Y} = [0,1]^{d_\rvy}$.
We have access to a dataset of $n$ observations $\mathcal{D}=\{(\vx_i, \vy_i)\}_{i\in[n]}$, where $[n]=\{1,\dots,n\}$,
which we split into disjoint subsets $\mathcal{D}=\mathcal{D}^{\textrm{train}}\cup\mathcal{D}^{\textrm{cal}}$ (no shuffling), of sizes $n^{\textrm{train}}$ and $n^{\textrm{cal}}$, respectively.
Throughout, we use $\hat{(\cdot)}$ to denote quantities learned or calibrated from data.

\paragraph{Score-based modeling.}
Let $s_\vtheta(\cdot): \mathcal{X}\times\mathbb{R}^{d_\rvy}\rightarrow \mathbb{R}_+$ be a partially convex function (convex in $\rvy$),
represented with neural networks (NNs) parameterized with weights $\vtheta$,
which we refer to as a \emph{score function}.
We use $s_\vtheta(\rvx,\rvy)$ to define an conditional density for $\rvy$ given $\rvx$,
\begin{equation}\label{eq:prob_score}
    p(\rvy|\rvx) = \frac{
    \text{exp}\big(-s_{\vtheta}(\rvx, \rvy)\big)
    }{Z_\vtheta(\rvx)} \propto \text{exp}\big(-s_{\vtheta}(\rvx, \rvy)\big),
\end{equation}
where $Z_\vtheta(\rvx)=\int_{\mathbb{R}^{d_{\rvy}}} \text{exp}\big(-s_{{\vtheta}}(\rvx, \tilde\rvy)\big)d\tilde\rvy$ is the partition function.
Given $\rvx$, a minimizer of $s_\vtheta(\rvx,\rvy)$ is a conditional mode under the score-induced density and provides a point forecast.
Learning consists of finding weights $\vtheta$ such that $s_\vtheta(\rvx,\rvy)$ associates low scores to correct values of $\rvy$ and high scores to incorrect ones \citep{lecun2006tutorial, song2021train}.

\paragraph{Prediction sets.}
We represent multivariate uncertainty using convex prediction sets defined as sublevel sets of the score function $s_\vtheta(\cdot)$, given by
\begin{equation}\label{set_form}
\mathcal{U}_{\vtheta}(\rvx;\rho)
=
\{\rvy\in\mathbb{R}^{d_{\rvy}}:
s_{\vtheta}(\rvx,\rvy)\leq\rho\}.
\end{equation}
The weights $\vtheta$ determine the conditional score model, while $\rho$ selects its sublevel set.
For fixed $\vtheta$ and $\rvx$, these sets are nested:
\begin{equation*}
\rho_a\leq\rho_b
\quad\Longrightarrow\quad
\mathcal{U}_{\vtheta}(\rvx;\rho_a)
\subseteq
\mathcal{U}_{\vtheta}(\rvx;\rho_b).
\end{equation*}
Thus, increasing $\rho$ enlarges the prediction set while retaining the trained score model.
Hereafter, we refer to $\rho$ as the \emph{score threshold}, or simply the threshold.

\paragraph{Robust optimization (RO).}
We consider parametric RO problems where the uncertainty $\rvy$ depends on $\rvx$ of the form
\begin{subequations}\label{rob_linear_prog}
\begin{align}
\underset{\vz}{\min}\,\,
& \vc^{\top}\vz, \\
\text{s.t.} \,\,
& g(\vz, \rvy) \leq 0, \quad \forall \rvy \in \mathcal{U}_\vtheta(\rvx; \rho),\label{robust_constr}
\end{align}
\end{subequations}
where $\vz$ denotes the decision variables, $\vc^\top \vz$ is a linear cost, the prediction set $\mathcal{U}_\vtheta(\rvx; \rho)$ coincides with an uncertainty set,
and $g(\vz,\rvy)$ is assumed to be the maximum of functions that are concave in $\rvy$ and convex in $\vz$.
Our motivating examples include robust linear inequalities of the form $g(\vz,\rvy)=\max_{l=1,\dots,L}\{a_l(\vz)^\top \rvy - b_l(\vz)\}$.
The robust constraint \eqref{robust_constr} is used to enforce {reliability requirements} \citep{bertsimas2021probabilistic}
by approximating the probabilistic chance constraint
\begin{equation}\label{chance_constr}
    \textrm{Prob}\big(g(\vz, \rvy) \leq 0\big) \geq 1-\varepsilon,
\end{equation}
where $\varepsilon>0$ is a (small) \emph{tolerance level} that controls the violation probability (e.g., $\varepsilon=5\%$).
Equivalently, we use the indicator function $\mathbb{I}\{\cdot\}$ to express \eqref{chance_constr} as
$\mathbb{E}[\mathbb{I}\{g(\vz,\rvy)\leq0\}]\geq1-\varepsilon$,
or, in terms of constraint violations, as
$\mathbb{E}[\mathbb{I}\{g(\vz,\rvy)>0\}]\leq\varepsilon$.
When $\mathcal{U}_\vtheta(\rvx; \rho)$ has a prescribed geometric form (e.g., a box, ellipsoid, or polyhedral set),
\eqref{robust_constr} can be reformulated into tractable deterministic constraints \citep{bertsimas_hertog2020}.

For a trained score model with weights $\hat\vtheta$, the threshold $\rho$ controls uncertainty-set size and, therefore, the cost--reliability trade-off induced by \eqref{rob_linear_prog}.
\emph{Coverage-based calibration} seeks the smallest $\rho$ such that
$\rvy \in \mathcal{U}_{\hat\vtheta}(\rvx;\rho)$ with probability at least $1-\varepsilon$,
which is directly analogous to
``maximizing the sharpness of predictive distributions subject to calibration''
\citep{gneiting2007probabilistic}.
Under robust feasibility, achieving $1-\varepsilon$ coverage is sufficient for \eqref{chance_constr}.
Coverage-based calibration can therefore be conservative when out-of-set realizations also satisfy the constraints.
In contrast, \emph{decision-based calibration} seeks the smallest $\rho$ that achieves the target reliability $1-\varepsilon$ for the resulting downstream decisions,
equivalently controlling the expected constraint-violation loss $\mathbb{E}\!\left[\mathbb{I}\{g(\vz^\star(\rvx;\rho),\rvy)>0\}\right]$
at level $\varepsilon$, where $\vz^\star(\rvx;\rho)$ denotes the solution of \eqref{rob_linear_prog}.

Although the prediction sets in \eqref{set_form} are nested in $\rho$, the downstream violation loss may not be monotone because the optimal decision $\vz^\star(\rvx;\rho)$ may change with the uncertainty set.
We address this issue through a monotonicity correction before applying CRC \citep{angelopoulos2024conformal}, as detailed in Section~\ref{sec:calibration}.

\subsection{Method Outline}\label{subsec:outline}

\begin{figure}[tb!]
    \centering

    \begin{minipage}[c]{0.49\columnwidth}
        \centering
        \resizebox{\linewidth}{!}{ \tikzset{every picture/.style={line width=0.75pt}} %set default line width to 0.75pt        

\begin{tikzpicture}[x=0.75pt,y=0.75pt,yscale=-1,xscale=1]
%uncomment if require: \path (0,295); %set diagram left start at 0, and has height of 295

%Straight Lines [id:da24431723046905673] 
\draw    (50,40) -- (117,40) ;
\draw [shift={(120,40)}, rotate = 180] [fill={rgb, 255:red, 0; green, 0; blue, 0 }  ][line width=0.08]  [draw opacity=0] (8.93,-4.29) -- (0,0) -- (8.93,4.29) -- cycle    ;
%Shape: Rectangle [id:dp3276230817017993] 
\draw   (127.5,20) -- (247.5,20) -- (247.5,60) -- (127.5,60) -- cycle ;

%Straight Lines [id:da10108306046134696] 
\draw    (260,40) -- (317,40) ;
\draw [shift={(320,40)}, rotate = 180] [fill={rgb, 255:red, 0; green, 0; blue, 0 }  ][line width=0.08]  [draw opacity=0] (8.93,-4.29) -- (0,0) -- (8.93,4.29) -- cycle    ;

%Shape: Rectangle [id:dp5061562558145523] 
\draw   (80,120) -- (295,120) -- (295,200) -- (80,200) -- cycle ;
%Straight Lines [id:da8755662568058671] 
\draw    (187.5,149.5) -- (187.5,166.5) ;
\draw [shift={(187.5,169.5)}, rotate = 270] [fill={rgb, 255:red, 0; green, 0; blue, 0 }  ][line width=0.08]  [draw opacity=0] (8.93,-4.29) -- (0,0) -- (8.93,4.29) -- cycle    ;

%Straight Lines [id:da35395287346826465] 
\draw    (40,160) -- (67,160) ;
\draw [shift={(70,160)}, rotate = 180] [fill={rgb, 255:red, 0; green, 0; blue, 0 }  ][line width=0.08]  [draw opacity=0] (8.93,-4.29) -- (0,0) -- (8.93,4.29) -- cycle    ;
%Straight Lines [id:da76824232510343] 
\draw    (300,160) -- (317,160) ;
\draw [shift={(320,160)}, rotate = 180] [fill={rgb, 255:red, 0; green, 0; blue, 0 }  ][line width=0.08]  [draw opacity=0] (8.93,-4.29) -- (0,0) -- (8.93,4.29) -- cycle    ;
%Right Arrow [id:dp42198518901497595] 
\draw   (192.88,209.74) -- (193.24,233.77) -- (198.97,233.68) -- (187.74,249.87) -- (176.03,234.02) -- (181.76,233.94) -- (181.41,209.91) -- cycle ;
%Straight Lines [id:da30624709088466484] 
\draw    (290,40) -- (290,80) -- (190,80) -- (190,107) ;
\draw [shift={(190,110)}, rotate = 270] [fill={rgb, 255:red, 0; green, 0; blue, 0 }  ][line width=0.08]  [draw opacity=0] (8.93,-4.29) -- (0,0) -- (8.93,4.29) -- cycle    ;

% Text Node
\draw (157,254) node [anchor=north west][inner sep=0.75pt]   [align=left] {$\displaystyle U_{\hat{\theta }}( x,\hat{\rho })$};
% Text Node
\draw (99,130.5) node [anchor=north west][inner sep=0.75pt]   [align=left] {Evaluate downstream loss};
% Text Node
\draw (83.5,171.5) node [anchor=north west][inner sep=0.75pt]   [align=left] {Monotonicity correction $\displaystyle +$ CRC};
% Text Node
\draw (140,29) node [anchor=north west][inner sep=0.75pt]   [align=left] {Train $\displaystyle s_{\hat{\mathbf{\theta }}}( x,y)$};
% Text Node
\draw (11,29.5) node [anchor=north west][inner sep=0.75pt]   [align=left] {$\displaystyle \mathcal{D}^{\text{train}}$};
% Text Node
\draw (320,31.5) node [anchor=north west][inner sep=0.75pt]   [align=left] {$\displaystyle \hat{\theta }$: location and shape};
% Text Node
\draw (11,149.5) node [anchor=north west][inner sep=0.75pt]   [align=left] {$\displaystyle \mathcal{D}^{\text{cal}}$};
% Text Node
\draw (320,150) node [anchor=north west][inner sep=0.75pt]   [align=left] {$\displaystyle \hat{\rho }$: set size};

\end{tikzpicture}
        }
        \smallskip
        \vspace{-.3in}
        \textbf{(a)}
    \end{minipage}
    \hfill
    \begin{minipage}[c]{0.49\columnwidth}
        \centering
        \includegraphics[width=\linewidth]{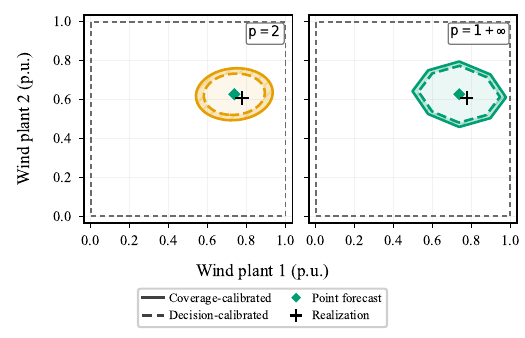}
        \smallskip
        \textbf{(b)}
    \end{minipage}
\caption{
Overview of the proposed methodology.
\textbf{(a)} Train the conditional score model using
$\mathcal{D}^{\textrm{train}}$ and calibrate the size threshold using
$\mathcal{D}^{\textrm{cal}}$.
\textbf{(b)} Coverage- and decision-calibrated prediction sets for two trained score models at $\varepsilon=5\%$.
The gray dashed boundary indicates the physical support
$\mathcal{Y}=[0,1]^2$.
}
    \label{flowchart}
\end{figure}

We prescribe the \emph{geometric form} of the prediction sets through the norm-based score specification.
Throughout, we use the term \emph{shape} to denote the orientation and context-dependent directional scaling,
whereas set location is specified by the conditional point predictor.
With the score model fixed, the threshold $\rho$ provides an additional scalar adjustment of set size.
We construct the prediction sets and use them in RO through two steps:
\begin{enumerate}
\item \textbf{Learn set location and shape.}
For the prescribed score form, we use
$\mathcal{D}^{\textrm{train}}$ to learn the weights $\vtheta$ that determine conditional set location and shape
(Section~\ref{sec:learning}).
\item \textbf{Calibrate the size threshold.}
With $\hat\vtheta$ fixed, we use
$\mathcal{D}^{\textrm{cal}}$ to evaluate downstream constraint
violations over a fixed threshold grid, apply the monotonicity
correction, and select $\rho$ using CRC at tolerance
$\varepsilon$
(Section~\ref{sec:calibration}).
\end{enumerate}

Figure~\ref{flowchart} summarizes the proposed approach.
Plot~(a) summarizes the learning and calibration steps.
Plot~(b) compares coverage- and decision-based calibration for two prescribed geometric forms, corresponding to the $2$-norm and $1+\infty$ scores.
For each trained score model, both procedures retain the same learned location and shape transformation,
with decision-based calibration selecting a smaller threshold and producing a smaller set in the illustrated cases.

\section{Learning Conditional Multivariate Prediction Sets}
\label{sec:learning}

In this section, we define the norm-based score functions
(Subsection~\ref{sec:score_design}) and describe how the corresponding score models are trained
(Subsection~\ref{sec:training}).

\subsection{Norm-based Score Design}
\label{sec:score_design}

Inspired by \cite{braun2025minimum}, we consider two NN-based predictive functions $\vf_\veta: \mathcal{X} \rightarrow \mathcal{Y}$ and $\mL_\vphi: \mathcal{X} \rightarrow \mathbb{R}^{d_{\rvy}\times d_{\rvy}}$ that are parameterized by weights $\veta$ and $\vphi$, respectively.
We define a $p$-norm score function
\begin{equation}\label{eq:norm_score}
s_{\vtheta}(\rvx, \rvy) \;=\; \left\Vert \mL_{\vphi}(\rvx)^{-1}\big(\rvy-\vf_{\veta}(\rvx)\big) \right\Vert_p,
\end{equation}
where $\vf_{\veta}(\rvx)$ is a conditional location predictor,
$\mL_{\vphi}(\rvx)$ is an invertible lower-triangular transformation with positive diagonal entries, $p\in\{1,2,\infty\}$ is the norm order, and
$\vtheta=(\veta,\vphi)$ collects the trainable weights.
The associated matrix
$\mSigma_{\vphi}(\rvx)=
\mL_{\vphi}(\rvx)\mL_{\vphi}(\rvx)^\top$
describes conditional dispersion.
The score function \eqref{eq:norm_score} measures the norm of a transformed residual
$\vu=\mL_\vphi(\rvx)^{-1}\big(\rvy-\vf_\veta(\rvx)\big)$,
where $\vf_\veta(\rvx)$ defines set location, while $\mL_\vphi(\rvx)$ determines set shape through a linear transformation of the prescribed norm ball.
\footnote{In power system applications where point predictors are provided by external vendors,
\eqref{eq:norm_score} or \eqref{eq:sum_scores_l1inf} can be used by fixing $\vf_{\veta}$ to the external predictor and learning only $\mL_{\vphi}$.}
Before intersection with the physical support, varying $\rho$
uniformly scales the resulting set about its learned location.

We focus on $p\in\{1,2,\infty\}$ because these choices yield convex prediction sets and tractable RO reformulations.
Before intersection with the physical support, the $2$-norm yields ellipsoidal sets, while the $1$- and $\infty$-norms yield polyhedral sets.

\paragraph{Sum of scores.}
To increase modeling flexibility while preserving convexity, we further employ a score function defined as the nonnegative sum of $K$ $p$-norm scores, given by
\begin{equation*}\label{eq:sum_scores_general}
s^{\textrm{sum}}_\vtheta(\rvx,\rvy)\;=\;\sum_{k\in[K]}\omega_k
\left\Vert \mL_\vphi(\rvx)^{-1}\big(\rvy-\vf_\veta(\rvx)\big) \right\Vert_{p_k},
\end{equation*}
where $\omega_k>0$ are fixed weights and the same $(\vf_{\veta},\mL_{\vphi})$ are shared across the terms.
In particular, we employ a score that equally weighs the $1$-norm and $\infty$-norm scores
\begin{align}\label{eq:sum_scores_l1inf}
s^{\textrm{sum}}_\vtheta(\rvx,\rvy)\;=\;
&
\left\Vert \mL_{\vphi}(\rvx)^{-1}\big(\rvy-\vf_{\veta}(\rvx)\big) \right\Vert_{1}
+ \nonumber\\
&
\left\Vert \mL_{\vphi}(\rvx)^{-1}\big(\rvy-\vf_{\veta}(\rvx)\big) \right\Vert_{\infty},
\end{align}
which creates polyhedral prediction sets that are more flexible than the ones created by either norm alone.
The $1$-norm component encourages the set to adapt to sparse deviations, where only a few variables have large forecast errors,
whereas the $\infty$-norm component controls the maximum componentwise deviation, avoiding sets that are too narrow along any individual dimension.
Their sum produces polyhedral prediction sets that can accommodate both concentrated and distributed error patterns -- see Fig.~\ref{flowchart}, plot (b) for an illustration.

{
The proposed family of scores \eqref{eq:norm_score}--\eqref{eq:sum_scores_l1inf} yields flexible convex prediction sets that capture multivariate dependence and admit tractable downstream RO reformulations, while supporting efficient likelihood-based training, which we analyze next.}

\subsection{Training Norm-based Score Models}
\label{sec:training}

We learn $\vtheta$ on $\mathcal{D}^{\textrm{train}}$ by applying gradient-based training and
minimizing a negative log-likelihood (NLL) objective induced by the score-based model.
For the $p$-norm score \eqref{eq:norm_score},
from \eqref{eq:prob_score}, the (unnormalized) conditional density has the form
$p_\vtheta(\rvy|\rvx)\propto \exp(-s_\vtheta\big(\rvx,\rvy)\big)$.
A change of variables $\vu=\mL_{\vphi}(\rvx)^{-1}\big(\rvy-\vf_{\veta}(\rvx)\big)$ yields a log-partition term of the form
$\log Z_{\vtheta}(\rvx)=\log|\det(\mL_{\vphi}(\rvx))| + \mathrm{const}(p,d_{\rvy})$,
where the constant term depends only on $(p,d_{\rvy})$ and does not affect optimization.

Dropping constants, for the $i$th observation, the NLL loss is given by
\begin{equation*}\label{eq:nll_pointwise}
\ell^{\textrm{NLL}}_{p,i}(\vtheta)
:=
s_{\vtheta}(\vx_i,\vy_i)
+
\log \big|\det(\mL_{\vphi}(\vx_i))\big|.
\end{equation*}
As $\mL_{\vphi}(\vx_i)$ is lower-triangular, $|\det(\mL_{\vphi}(\vx_i))|$ is efficiently computed as the product of its diagonal entries.
Similarly, the NLL for the sum score \eqref{eq:sum_scores_l1inf} is given by
\begin{equation*}\label{eq:nll_pointwise_sum}
\ell^{\textrm{NLL}}_{\textrm{sum},i}(\vtheta)
:=
s^{\textrm{sum}}_{\vtheta}(\vx_i,\vy_i)
+
\log \big|\det(\mL_{\vphi}(\vx_i))\big|.
\end{equation*}
When $\vf_\veta$ and $\mL_\vphi$ are trained jointly, we also include the mean squared error to stabilize training and improve point forecast accuracy, given by
\begin{equation*}\label{eq:mse_pointwise}
\ell^{\textrm{MSE}}_{i}(\veta)
:=
\left\Vert \vy_i - \vf_{\veta}(\vx_i)\right\Vert_2^2,
\end{equation*}
which is weighted by a user-defined hyperparameter $\gamma \geq 0$.
For the $p$-norm score \eqref{eq:norm_score}, the final loss is given by
\begin{equation}\label{eq:final_train_obj}
 \ell_{p,i}^{\textrm{NLL}}(\vtheta) + \gamma \ell_{i}^{\textrm{MSE}}(\veta),
\end{equation}
and is derived similarly for the sum score \eqref{eq:sum_scores_l1inf} by replacing $\ell_{p,i}^{\textrm{NLL}}(\vtheta)$ with $\ell_{\textrm{sum},i}^{\textrm{NLL}}(\vtheta)$.

To ensure that $\mL_{\vphi}(\rvx)$ is invertible, we use a lower-triangular parameterization and pass its diagonal outputs through the softplus activation $\textrm{softplus}(u)=\ln(1+\exp(u))$ \citep{paszke2019pytorch}.
In addition, for $p\in\{1,\infty\}$, the score function is non-differentiable.
During gradient-based training, we use a smooth approximation (e.g., log-sum-exp smoothing for $\|\cdot\|_\infty$ and a smooth absolute value for $\|\cdot\|_1$), and evaluate the exact score at inference.
During training, we optimize \eqref{eq:final_train_obj} without imposing the support constraints $\mathcal{Y}=[0,1]^{d_{\rvy}}$ on $\rvy$,
which allows us to retain the closed-form partition function, and hence the tractable NLL objective derived above.
As the score-based model defines an unnormalized density proportional to $\exp(-s_\vtheta(\rvx,\rvy))$, the learned density is not explicitly restricted to the physical support.
To ensure feasible prediction sets, we apply a two-step post-processing at inference,
where we first project the output of $\vf_{\hat \veta}$ onto $[0,1]^{d_{\rvy}}$ and then intersect the prediction set with the known support;
hereafter, $\vf_{\hat\veta}$ denotes the projected point predictor.
The final \emph{feasible} prediction set is given by
\begin{equation} \label{feasible_set_form}
    \mathcal{U}_{\hat\vtheta}^{\textrm{feas}}(\rvx; \rho) = \mathcal{U}_{\hat\vtheta}(\rvx; \rho)\cap [0,1]^{d_{\rvy}}.
\end{equation}

\section{Robust DCOPF for Reserve Deliverability}
\label{sec:dcopf}

In this section, we motivate and formulate the reserve deliverability problem (Subsection~\ref{sec:dcopf_formulation}),
and provide tractable robust reformulations (Subsection~\ref{sec:dcopf_rc}).

\subsection{Problem Setting and Formulation}\label{sec:dcopf_formulation}

In power system operations, SOs schedule generation ahead of delivery through market clearing and procure reserves to maintain the real-time balance between supply and demand \citep{ferc2021}.
However, network congestion may prevent scheduled reserves from balancing realized forecast errors, even when sufficient aggregate reserve capacity is available \citep{caiso2020deliverability, chen2023addressing}, increasing costs and threatening system reliability.
Here, reserve deliverability means that the scheduled recourse policy can balance a realized forecast error without exceeding procured reserves or transmission limits.
Ensuring deliverability therefore requires jointly considering energy, reserves, and transmission capacity \citep{vancaelenberg2026deployment, madani2026reservedeliverability}.
We consider an idealized market setting in which energy and reserves are scheduled 15 minutes ahead using a DCOPF approximation.
The ARO formulation jointly determines the schedule and an affine recourse policy that accommodates forecast errors within the uncertainty set \citep{vancaelenberg2026deployment}.
Figure~\ref{fig:operating_timeline} illustrates the operational plan,
where features enter the trained score model to construct a prediction set, which then informs the energy and reserve schedule and the affine recourse policy.
When renewable production is realized, this policy determines reserve activation in response to forecast errors.
New observations update the features for the next scheduling cycle, repeated every 15 minutes.

\begin{figure}[tb!]
    \centering
    \resizebox{0.5\columnwidth}{!}{%
   \tikzset{every picture/.style={line width=0.75pt}} %set default line width to 0.75pt        

\begin{tikzpicture}[x=0.75pt,y=0.75pt,yscale=-1,xscale=1]
%uncomment if require: \path (0,277); %set diagram left start at 0, and has height of 277

%Straight Lines [id:da3004804715338182] 
\draw    (500,200) -- (500,220) -- (80,219) -- (80,193) ;
\draw [shift={(80,190)}, rotate = 90] [fill={rgb, 255:red, 0; green, 0; blue, 0 }  ][line width=0.08]  [draw opacity=0] (8.93,-4.29) -- (0,0) -- (8.93,4.29) -- cycle    ;
%Straight Lines [id:da03374484085560148] 
\draw    (82.5,121) -- (82.5,148) ;
\draw [shift={(82.5,151)}, rotate = 270] [fill={rgb, 255:red, 0; green, 0; blue, 0 }  ][line width=0.08]  [draw opacity=0] (8.93,-4.29) -- (0,0) -- (8.93,4.29) -- cycle    ;

%Straight Lines [id:da12083091005272806] 
\draw    (147.63,165.25) -- (254.63,165.25) ;
\draw [shift={(257.63,165.25)}, rotate = 180] [fill={rgb, 255:red, 0; green, 0; blue, 0 }  ][line width=0.08]  [draw opacity=0] (8.93,-4.29) -- (0,0) -- (8.93,4.29) -- cycle    ;
%Straight Lines [id:da6397091806134982] 
\draw    (339.89,165.25) -- (446.89,165.25) ;
\draw [shift={(449.89,165.25)}, rotate = 180] [fill={rgb, 255:red, 0; green, 0; blue, 0 }  ][line width=0.08]  [draw opacity=0] (8.93,-4.29) -- (0,0) -- (8.93,4.29) -- cycle    ;
%Straight Lines [id:da15408109703661033] 
\draw    (484.5,112) -- (484.5,139) ;
\draw [shift={(484.5,142)}, rotate = 270] [fill={rgb, 255:red, 0; green, 0; blue, 0 }  ][line width=0.08]  [draw opacity=0] (8.93,-4.29) -- (0,0) -- (8.93,4.29) -- cycle    ;
%Shape: Rectangle [id:dp8598636479069699] 
\draw   (264.26,140.5) -- (333.26,140.5) -- (333.26,190) -- (264.26,190) -- cycle ;

%Shape: Rectangle [id:dp039786757342001766] 
\draw   (450,142) -- (530,142) -- (530,191.5) -- (450,191.5) -- cycle ;

% Text Node
\draw (147,231) node [anchor=north west][inner sep=0.75pt]   [align=left] {New observations; repeat every 15 min};
% Text Node
\draw (154,134) node [anchor=north west][inner sep=0.75pt]   [align=left] {Set $\displaystyle \mathcal{U}_{\hat{\vtheta }}( \vx_{0} ;\hat{\rho})$};
% Text Node
\draw (457.5,146.25) node [anchor=north west][inner sep=0.75pt]   [align=left] {Reserve \\activation};
% Text Node
\draw    (27,150.25) -- (143,150.25) -- (143,180.25) -- (27,180.25) -- cycle  ;
\draw (30,154.25) node [anchor=north west][inner sep=0.75pt]   [align=left] {Model $\displaystyle s_{\hat{\vtheta}}( \vx_{0} ,\rvy)$ };
% Text Node
\draw (348,122) node [anchor=north west][inner sep=0.75pt]   [align=left] {Schedule, \\linear policy};
% Text Node
\draw (141,10) node [anchor=north west][inner sep=0.75pt]   [align=left] {Scheduling at time $\displaystyle t$};
% Text Node
\draw (391,10) node [anchor=north west][inner sep=0.75pt]   [align=left] {Delivery at time $\displaystyle t+15$};
% Text Node
\draw (42,101) node [anchor=north west][inner sep=0.75pt]   [align=left] {Features $\displaystyle \vx_{0}$};
% Text Node
\draw (441,62) node [anchor=north west][inner sep=0.75pt]   [align=left] {Realized error\\ $\displaystyle \vxi =\vy-\hat{\vy}$};
% Text Node
\draw (271.76,146.25) node [anchor=north west][inner sep=0.75pt]   [align=left] {Robust \\DCOPF};

\end{tikzpicture}
    }
    \caption{
    Flowchart of the operating timeline, including forecasting, robust scheduling, and reserve activation at 15-min intervals.
    }
    \label{fig:operating_timeline}
\end{figure}

Consider a power system where $\mathcal{N}$ is the set of nodes, $\mathcal{L}$ is the set of lines, $\mathcal{G}$ is the set of dispatchable generators, and $\mathcal{R}$ is the set of stochastic renewable generators.
We assume that the nodal demand $\vd\in\mathbb{R}^{|\mathcal{N}|}$ is a deterministic parameter known at the scheduling stage and consider uncertainty only in renewable production.
The framework extends directly to uncertain demand by representing uncertainty as nodal net load (demand minus renewable production).
Let $\rvy = \hat\vy + \vxi$ be the stochastic renewable production,
where $\hat\vy = \vf_{\hat\veta}(\rvx) \in\mathbb{R}^{|\mathcal{R}|}$ denotes the point forecast,
and $\vxi\in\mathbb{R}^{|\mathcal{R}|}$ is the random forecast error.
Given features $\rvx$ and fixed score weights $\hat\vtheta$, our method outputs a prediction set for renewable production
$\mathcal{U}^{\textrm{feas}}(\rvx;\rho)$ of the form \eqref{feasible_set_form} (the dependency on $\hat\vtheta$ is suppressed for simplicity).
For convenience, we express uncertainty in terms of forecast errors $\vxi$ centered on the point forecast, using the map
\begin{align}
    \label{eq:V_norm}
\mathcal{V}(\rvx;\rho) = \left\{\vxi:
\ \vf_{\hat\veta}(\rvx) + \vxi \in \mathcal{U}^{\textrm{feas}}(\rvx;\rho) \right\}.
\end{align}
For instance, for the $p$-norm score \eqref{eq:norm_score}, the feasible forecast-error prediction set is given by
$\mathcal{V}(\rvx;\rho) = \left\{\vxi:\ \left\|\mL_{\hat \vphi}(\rvx)^{-1}\vxi\right\|_{p}\le \rho, \; -\vf_{\hat\veta}(\rvx)\le \vxi \le \vone-\vf_{\hat\veta}(\rvx)\right\}$.
The process for the sum score \eqref{eq:sum_scores_l1inf} is similar.
\footnote{The score models are trained using renewable production normalized by nominal plant capacities.
The forecast errors are rescaled to nominal units before entering the constraints; we absorb this inverse scaling into $\mathcal V(\rvx;\rho)$ to keep the notation compact.}

Following \citet{mieth2023prescribed}, we approximate the two-stage ARO problem \citep{vancaelenberg2026deployment} using linear decision rules, resembling automatic generation control \citep{hou2020chance}.
Let $\vp\in\mathbb{R}^{|\mathcal{G}|}$ denote the dispatch schedule and
$\mA\in\mathbb{R}^{|\mathcal{G}|\times|\mathcal{R}|}$ the recourse matrix.
For an error realization $\vxi$, generator outputs follow the affine policy
\begin{equation*}
\vp(\vxi)=\vp-\mA\vxi.
\end{equation*}
The robust DCOPF jointly determines the energy schedule, reserve capacities,
and recourse matrix under linear costs:
\begin{subequations}\label{dcopf}
\begin{align}
\underset{\substack{\vp, \vr^+, \vr^-, \\\vm^{+}, \vm^{-},
\mA}}{\min} & \,\, (\vc^{\textrm{e}})^{\top}\vp  + (\vc^{\textrm{r}})^{\top}(\vr^{+} + \vr^{-}), \label{dsw_obj} \\
\text{s.t.} & \, \, \vone^{\top}\vp + \vone^{\top}\hat\vy = \vone^{\top}\vd,  \label{dcopf_c1} \\
& \, \,  \mA^{\top}\vone = \vone, \label{dcopf_c2} \\
& \, \, \mB^{\mathcal{G}} \vp +\mB^{\mathcal{R}} \hat\vy - \mB^{\mathcal{N}}\vd = \vf^{\max}-\vm^{+}, \label{dcopf_c3}\\
 & \, \, -\big(\mB^{\mathcal{G}} \vp +\mB^{\mathcal{R}} \hat\vy - {\mB^{\mathcal{N}}}\vd\big) = \vf^{\max}-\vm^{-},\label{dcopf_c4}\\
 & \, \, \vp + \vr^+ \leq \vp^{\max}, \label{dcopf_c5}\\
& \,\, \vp - \vr^- \geq \vp^{\min},  \label{dcopf_c6}\\
& \, \, \vp, \vr^+, \vr^-, \vm^+, \vm^- \geq \vzero,  \label{dcopf_c7}
\\
& -\mA \vxi \leq \vr^+,  \hspace{4.85em}{\quad \forall\, \vxi \in \mathcal{V}(\rvx;\rho)}, \label{dcopf_c8}\\
& \mA \vxi \leq \vr^-,
\hspace{5.75em}{\quad \forall\, \vxi \in  \mathcal{V}(\rvx;\rho)}, \label{dcopf_c9}
\\
& (\mB^{\mathcal{R}} - \mB^{\mathcal{G}}\mA) \vxi \leq \vm^+,
\hspace{1em}{\quad \forall\, \vxi \in \mathcal{V}(\rvx;\rho)}, \label{dcopf_c10}
\\
& -(\mB^{\mathcal{R}} - \mB^{\mathcal{G}}\mA) \vxi \leq \vm^-,
{\quad \forall\, \vxi \in \mathcal{V}(\rvx;\rho)}.
\label{dcopf_c11}
\end{align}
\end{subequations}
The objective \eqref{dsw_obj} minimizes energy and reserve procurement costs,
where $\vc^{\textrm{e}}$ and $\vc^{\textrm{r}}$ are the respective cost coefficients and $\vr^+,\vr^-\in\mathbb{R}^{|\mathcal{G}|}$ denote upward and downward reserve capacities.
Constraint \eqref{dcopf_c1} balances scheduled generation and demand,
while \eqref{dcopf_c2} ensures that total recourse offsets the renewable forecast error.
Constraints \eqref{dcopf_c3}--\eqref{dcopf_c4} define the non-negative transmission margins $\vm^+,\vm^-\in\mathbb{R}^{|\mathcal{L}|}$ around the scheduled flows, and \eqref{dcopf_c5}--\eqref{dcopf_c6} keep energy and reserve schedules within generator limits.
Here, $\vp^{\min}$, $\vp^{\max}$, and $\vf^{\max}$ denote generation and transmission limits,
while $\mB^{\mathcal{G}}$, $\mB^{\mathcal{R}}$, and $\mB^{\mathcal{N}}$
combine power transfer distribution factors with the corresponding node incidence matrices.
Finally, the robust constraints \eqref{dcopf_c8}--\eqref{dcopf_c9} keep generator adjustments within procured reserves,
and \eqref{dcopf_c10}--\eqref{dcopf_c11} keep the resulting flow changes within available transmission margins,
for all forecast errors in $\mathcal{V}(\rvx;\rho)$.

\subsection{Robust Constraint Reformulation}
\label{sec:dcopf_rc}

We now show how to reformulate the robust constraints using support functions -- see \citep[Chapter~2]{bertsimas_hertog2020} for background.
For notational simplicity, we write $\mL=\mL_{\hat\vphi}(\rvx)$ and suppress the dependence of the sets on $\rvx$.
Each robust constraint in \eqref{dcopf_c8}-\eqref{dcopf_c11} has the generic form
\begin{equation}\label{eq:generic_robust_linear}
\va^\top \vxi \le b,
\quad \forall \vxi\in \mathcal{S}(\rho),
\end{equation}
where both $\va, b$ depend on the decision variables (rows of $\mA$ or $(\mB^{\mathcal{R}}-\mB^{\mathcal{G}}\mA)$ and $\vr^\pm,\vm^\pm$, respectively) and $\mathcal{S}(\rho)$ is a convex set whose size depends on $\rho$.
Equivalently,
\begin{equation}\label{eq:support_fn_form}
\delta^*\big(\va\mid \mathcal{S}(\rho)\big) \le b,
\end{equation}
where $\delta^*\big(\va\mid \mathcal{S}(\rho)\big) := \sup_{\vxi\in\mathcal{S}(\rho)}\va^\top\vxi$ is the support function of $\mathcal{S}(\rho)$.
Hence, reformulating the robust constraints reduces to evaluating this support function.

\paragraph{$p$-norm prediction sets.}
If $\mathcal{S}^{}(\rho)=\{\vxi:\|\mL^{-1}\vxi\|_p\le \rho\}$ (no box truncation), then
\begin{equation}\label{eq:support_norm_only}
\sup_{\|\mL^{-1}\vxi\|_p\le \rho}\ \va^\top \vxi
=
\rho\,\big\|\mL^{\top}\va\big\|_{p^\ast},
\end{equation}
where $p^\ast$ is the dual norm exponent ($1/p+1/p^\ast=1$).
Substituting \eqref{eq:support_norm_only} into \eqref{eq:support_fn_form} replaces the robust constraint \eqref{eq:generic_robust_linear} with a finite constraint, which, depending on the choice of $p$, can be represented using standard linear or second-order cone reformulations.

\paragraph{Box prediction sets.}
Suppose now that the prediction set is an axis-aligned box, $\mathcal{S}=\{\vxi:\vl\le \vxi \le \vh\}$,
where $\vl$ and $\vh$ denote componentwise lower and upper bounds.
Since the maximization separates across coordinates,
\begin{align}
\sup_{\vl\le \vxi\le \vh} \va^\top \vxi & = \sum_{j=1}^{|\mathcal{R}|}\sup_{l_j\le \xi_j\le h_j} a_j \xi_j \nonumber\\
&= \sum_{j=1}^{|\mathcal{R}|} \max\{a_j h_j,\; a_j l_j\} \nonumber\\
&= \sum_{j=1}^{|\mathcal{R}|}
\big(\max\{a_j,0\}h_j + \min\{a_j,0\}l_j\big).
\label{eq:support_box}
\end{align}
Substituting \eqref{eq:support_box} into \eqref{eq:support_fn_form} replaces the robust constraint \eqref{eq:generic_robust_linear} with a linear programming reformulation,
using auxiliary variables to represent the maximum terms.

\paragraph{$p$-norm and box intersection.}
For the intersection set $\mathcal{S}(\rho)=\mathcal{S}_1(\rho)\cap \mathcal{S}_2$ of the form \eqref{eq:V_norm},
where $\mathcal{S}_1(\rho)=\{\vxi:\|\mL^{-1}\vxi\|_p\le \rho\}$ and $\mathcal{S}_2=\{\vxi:\vl \le \vxi\le \vh\}$,
the support function can be expressed as
\begin{align}\label{eq:support_intersection}
\delta^*(\va\mid \mathcal{S}(\rho))
=
\min_{\vw}\ & \delta^*(\vw\mid \mathcal{S}_1(\rho)) \nonumber\\ & +\delta^*(\va-\vw\mid \mathcal{S}_2),
\end{align}
using \citet[Eq. (2.38)]{bertsimas_hertog2020}.
Introducing $\vw$ as an auxiliary variable yields a finite robust counterpart by combining \eqref{eq:support_norm_only} and \eqref{eq:support_box} and noting that the $\min$ operator can be omitted, as we have a $\leq$ inequality.

\paragraph{Sum of $p$-norm scores.}
If the prediction set is induced by the sum of the $1$-norm and $\infty$-norm scores \eqref{eq:sum_scores_l1inf},
then we derive the support function of the polyhedral set
$\{\vxi:\|\mL^{-1}\vxi\|_1+\|\mL^{-1}\vxi\|_\infty\le \rho\}$ using duality.
The reformulated robust constraint remains tractable (polyhedral) and can be combined with the same intersection-based reformulation in \eqref{eq:support_intersection}.

\section{Decision-based Calibration for Score Threshold}
\label{sec:calibration}

This section describes how to calibrate the score threshold $\rho$, which determines the prediction-set size.
We present the CRC framework (Subsection~\ref{sec:calibration_toolkit}),
formulate calibration losses for the robust DCOPF problem (Subsection~\ref{sec:dcopf_losses}),
and develop the proposed calibration algorithm (Subsection~\ref{sec:rho_tuning}).

\subsection{Conformal Risk Control}
\label{sec:calibration_toolkit}

We first present CRC for a generic calibration loss $\ell$;
Subsection~\ref{sec:dcopf_losses} specializes this loss to predictive miscoverage and operational constraint violations.
We formulate calibration in terms of the forecast-error set $\mathcal V(\rvx;\rho)$, obtained from $\mathcal U^{\textrm{feas}}(\rvx;\rho)$ by \eqref{eq:V_norm}.
For $\rho\in[0,\rho^{\max}]$, let
$\ell(\mathcal{V}(\rvx;\rho),\vxi)\in[0,\beta]$.
For observation $(\vx_i,\vy_i)$, define the forecast error
$\vxi_i:=\vy_i-\vf_{\hat\veta}(\vx_i)$ and the loss function
$\ell_i(\rho):=
\ell(\mathcal{V}(\vx_i;\rho),\vxi_i)$.
CRC assumes that these loss functions are non-increasing and right-continuous in $\rho$.
For a fixed threshold, the population risk is the expected loss
\begin{equation*}
R(\rho):=
\mathbb{E}\!\left[
\ell(\mathcal{V}(\rvx;\rho),\vxi)
\right],
\end{equation*}
conditional on the trained model.
The empirical calibration risk estimates this expectation using the held-out observations:
\begin{equation}\label{eq:crc_defs}
\hat R(\rho):=
\frac{1}{n^{\textrm{cal}}}
\sum_{i\in[n^{\textrm{cal}}]}\ell_i(\rho).
\end{equation}

For a risk tolerance $\alpha\in(0,\beta)$, assume that
$\ell_i(\rho^{\max})\leq\alpha$ almost surely, including for a new test observation.
This terminal-loss condition ensures that the largest threshold achieves the prescribed tolerance;
for a binary loss and $\alpha<1$, it requires zero terminal loss almost surely.
CRC selects the smallest threshold satisfying the following
finite-sample correction to empirical risk:
\begin{equation}\label{eq:crc_rho_hat}
\hat\rho
\;:=\;
\inf\left\{
\rho:\ \frac{n^{\textrm{cal}}}{n^{\textrm{cal}}+1} \hat R(\rho)+\frac{\beta}{n^{\textrm{cal}}+1}\le \alpha
\right\},
\end{equation}
with $\hat\rho=\rho^{\max}$ if the feasible set is empty.
Under data exchangeability \citep[Theorem~1]{angelopoulos2024conformal}, this rule guarantees
\begin{equation}\label{eq:crc_guarantee}
\mathbb{E}\big[\ell_{n^{\textrm{cal}}+1}(\hat\rho)\big]\ \le\ \alpha.
\end{equation}
The expectation is jointly over the calibration observations and a new test observation, conditional on the data used to train the model.
Data exchangeability requires their conditional joint distribution to be invariant under data ordering.
For binary losses, $\beta=1$, so $\alpha$ can be interpreted as the tolerated probability of the event measured by the loss, such as downstream constraint violation.
The additive term $\beta/(n^{\textrm{cal}}+1)$ vanishes as $n^{\textrm{cal}}$ grows.

For a bounded loss that is \emph{non}-monotone in $\rho$,
\citet{angelopoulos2024conformal} further propose replacing the original loss with its monotone upper envelope.
Namely, for a loss $\ell_i(\rho)$, define
\begin{equation}\label{eq:monotone_majorant_general}
    \ell_i^{\uparrow}(\rho)
    :=
    \sup_{\rho'\geq \rho}\ell_i(\rho').
\end{equation}
By construction, $\ell_i^{\uparrow}(\rho)$ is monotone non-increasing and satisfies $\ell_i(\rho)\leq \ell_i^{\uparrow}(\rho)$;
applying CRC to $\ell_i^{\uparrow}$ therefore controls the expected value of the original loss as well.
We specialize this construction to the downstream constraint violation loss in Subsection~\ref{sec:rho_tuning}.

\subsection{Calibration Losses for Robust DCOPF}
\label{sec:dcopf_losses}

We now define calibration losses for the robust DCOPF problem \eqref{dcopf}, which can also be extended to other robust optimization problems.
Let $\vz = (\vp,\vr^{+},\vr^{-},\vm^{+},\vm^{-},\mA)$ denote the decision variables and define the maximum constraint residual
associated with the robust constraints \eqref{dcopf_c8}--\eqref{dcopf_c11}:
\begin{align}\label{eq:gDC_def}
g^{\textrm{DC}}(\vz, \vxi) :=
& \max\Big\{
\max_{g\in\mathcal{G}}\big[(-(\mA\vxi)_g- r_{g}^{+})\big],
\nonumber\\
&
\quad\max_{g\in\mathcal{G}}\big[((\mA\vxi)_g- r_{g}^{-})\big],\nonumber\\
&
\quad\max_{l\in\mathcal{L}}\big[\big((\mB^{\mathcal{R}}-\mB^{\mathcal{G}}\mA)\vxi\big)_l-m_{l}^{+}\big], \nonumber\\
&
\quad\max_{l\in\mathcal{L}}\big[\big(-(\mB^{\mathcal{R}}-\mB^{\mathcal{G}}\mA)\vxi\big)_l-m_{l}^{-}\big]
\Big\}.
\end{align}
Then \eqref{dcopf_c8}--\eqref{dcopf_c11} are equivalent to $g^{\textrm{DC}}(\vz,\vxi)\le 0$ for all $\vxi\in\mathcal{V}(\rvx;\rho)$.
For the thresholds considered below, assume that \eqref{dcopf} admits an optimal solution satisfying the original robust constraints with probability one and
let $\vz^\star(\rvx;\rho)\in \arg\min \eqref{dcopf}$ denote this solution.

\paragraph{Constraint violation loss.}
The decision-based calibration loss is the realized constraint violation indicator
\begin{equation}\label{eq:violation_indicator}
\ell^{\mathrm{viol}}_i(\rho)
:=
\mathbb{I}\Big\{ g^{\textrm{DC}}\big(\vz^\star(\vx_i;\rho),\vxi_i\big) > 0
\Big\}.
\end{equation}
Each observation contributes one violation if any reserve or transmission constraint is violated, irrespective of whether its forecast error lies inside the prediction set.
Since the expectation of this indicator equals the corresponding constraint violation probability, controlling
$\mathbb{E}[\ell_i^{\mathrm{viol}}(\rho)]$ at level $\varepsilon$ is equivalent to targeting the probabilistic constraint \eqref{chance_constr}.
However, $\ell_i^{\mathrm{viol}}(\rho)$ is generally not monotone in $\rho$.
As discussed, although increasing $\rho$ enlarges the uncertainty set $\mathcal{V}(\vx_i;\rho)$,
it also changes the optimal schedule and recourse policy, so the realized violation loss may not decrease monotonically.
We therefore apply CRC to a monotone upper envelope, detailed in Subsection~\ref{sec:rho_tuning}.

\paragraph{Miscoverage loss.}
We also consider the standard miscoverage loss
\begin{align}\label{eq:miscoverage_loss}
\ell^{\mathrm{mis}}_i(\rho)
& =
\mathbb{I}\{\vxi_i \notin \mathcal{V}(\vx_i;\rho)\} \nonumber \\
&=
\mathbb{I}\{\vy_i \notin \mathcal{U}^{\textrm{feas}}(\vx_i;\rho)\} =
\mathbb{I}\{s_{\hat\vtheta}(\vx_i, \vy_i) > \rho \},
\end{align}
where the last equality uses the fact that $\vy_i$ lies within the physical support.
The loss is non-increasing in $\rho$.
For tolerance $\varepsilon$, CRC recovers split CP: $\hat\rho^{\textrm{mis}}$ is the $k$th smallest calibration score, with $k={\lceil (n^{\textrm{cal}}+1)(1-\varepsilon)\rceil}$.
The resulting prediction set satisfies
\begin{equation}\label{eq:coverage_def}
\mathrm{Prob}\big(\vxi \in \mathcal{V}(\rvx;\hat\rho^{\textrm{mis}})\big) \ge 1-\varepsilon.
\end{equation}

The miscoverage loss also upper-bounds the constraint violation loss pointwise.
If $\vxi_i\in \mathcal{V}(\vx_i;\rho)$, then by robust feasibility of $\vz^\star(\vx_i;\rho)$ we have $g^{\mathrm{DC}}\bigl(\vz^\star(\vx_i;\rho),\vxi_i\bigr)\le 0$, and therefore $\ell_i^{\mathrm{viol}}(\rho)\le \ell_i^{\mathrm{mis}}(\rho)$ for any $\rho$.
Moreover, the coverage guarantee \eqref{eq:coverage_def} implies
\begin{equation*}\label{eq:coverage_implies_cc}
\begin{aligned}
&\mathrm{Prob}\!\left(
g^{\mathrm{DC}}\!\left(
\vz^\star(\rvx;\hat\rho^{\textrm{mis}}),\vxi
\right)\leq0
\right)
\\
&\qquad\geq
\mathrm{Prob}\!\left(
\vxi\in\mathcal{V}(\rvx;\hat\rho^{\textrm{mis}})
\right)
\geq1-\varepsilon.
\end{aligned}
\end{equation*}
See \citet[Lemma~1]{hong2021learning}.
Hence, $\hat\rho^{\textrm{mis}}$ provides a conservative coverage-based benchmark for the downstream reliability target.

\subsection{Decision-based Calibration with Monotonicity Correction}
\label{sec:rho_tuning}

We construct a monotone upper envelope by assigning a violation at a threshold whenever the same observation produces a violation at that threshold or \emph{any} larger candidate.
We use a fixed, ordered grid
\begin{equation*}
\mathcal{T} =
\{\rho^1<\rho^2<\cdots<\rho^M\},
\end{equation*}
where we assume that the largest candidate threshold $\rho^M$ is large enough such that $\mathcal{U}_{\hat\vtheta}^{\textrm{feas}}(\rvx;\rho^M) =[0,1]^{d_{\rvy}}$ with probability one.
That is, the terminal prediction set recovers the full physical support of renewable production.
Under the robust-feasibility assumption above, this implies
\begin{equation}\label{eq:crc_endpoint}
\ell^{\textrm{viol}}(\rho^M)=0,
\end{equation}
which is the CRC endpoint condition for the binary violation loss.
Both the threshold grid $\mathcal{T}$ and the trained score model $s_{\hat\vtheta}$ are fixed before evaluating the calibration losses.
In the numerical experiments, we select $\rho^M$ using training
and validation data before evaluating calibration losses, as described in Section~\ref{sec:exp_setup}.
The finite-sample guarantee below
requires the resulting grid to satisfy the endpoint condition
\eqref{eq:crc_endpoint}.

For observation $i$ and candidate threshold $\rho^j$ from the grid $\mathcal{T}$,
we first evaluate the raw constraint violation loss $\ell_i^{\textrm{viol}}(\rho^j)$ \eqref{eq:violation_indicator} and define the monotonicity correction as the following upper-envelope transformation
\begin{equation}
\label{eq:corrected_violation_loss}
\ell_i^{\textrm{viol},\uparrow}(\rho^j) :=
\max_{k=j,\ldots,M} \ell_i^{\textrm{viol}}(\rho^k).
\end{equation}
For example, if $\ell_i^{\textrm{viol}}(\rho^j)$ returns the sequence $(0,1,0)$ for $\rho^1<\rho^2<\rho^3$, then $\ell_i^{\textrm{viol}, \uparrow}(\rho^j)$ returns $(1,1,0)$;
similarly, an already non-increasing sequence remains unchanged.
The corrected loss $\ell_i^{\textrm{viol}, \uparrow}(\rho^j)$ is the smallest monotone non-increasing function on the grid that upper-bounds the raw loss:
\begin{equation}
\label{eq:viol_majorization}
\ell_i^{\textrm{viol}}(\rho^j)
\leq
\ell_i^{\textrm{viol},\uparrow}(\rho^j),
\qquad j=1,\ldots,M.
\end{equation}

Moreover,
the corrected violation loss \eqref{eq:corrected_violation_loss} remains pointwise bounded by the miscoverage loss:
\begin{equation}
\label{eq:loss_ordering}
\ell_i^{\textrm{viol}}(\rho^j)
\leq
\ell_i^{\textrm{viol},\uparrow}(\rho^j)
\leq
\ell_i^{\textrm{mis}}(\rho^j).
\end{equation}
Indeed, for every $k\geq j$, robust feasibility and nestedness give
$\ell_i^{\textrm{viol}}(\rho^k)\leq\ell_i^{\textrm{mis}}(\rho^k)\leq\ell_i^{\textrm{mis}}(\rho^j)$.
Taking the maximum over $k \geq j$ yields \eqref{eq:loss_ordering}.
Thus, the monotonicity correction can increase the violation loss but cannot exceed miscoverage at the same threshold.
This construction follows the monotonization procedure proposed by \citet{angelopoulos2024conformal} for general non-monotone losses.

We next define the corrected empirical risk as
\begin{equation}
\label{eq:majorized_risk}
\hat R^\uparrow(\rho^j)
=
\frac{1}{n^{\textrm{cal}}}
\sum_{i\in[n^{\textrm{cal}}]}
\ell_i^{\textrm{viol},\uparrow}(\rho^j).
\end{equation}
Since $\ell_i^{\textrm{viol},\uparrow}(\rho^j)\in\{0,1\}$, we have $\beta=1$ in \eqref{eq:crc_rho_hat}.
Therefore, for the target violation tolerance $\varepsilon$, the CRC criterion reduces to
\begin{equation}
\label{eq:grid_crc_condition}
\hat R^\uparrow(\rho^j)
\leq
\varepsilon-\frac{1-\varepsilon}{n^{\textrm{cal}}}.
\end{equation}
We select the smallest candidate threshold satisfying \eqref{eq:grid_crc_condition}.
Under \eqref{eq:crc_endpoint}, $\rho^M$ satisfies this criterion whenever $\varepsilon\geq1/(n^{\textrm{cal}}+1)$, so the feasible candidate set is non-empty regardless of grid spacing.

\begin{proposition}[Finite-sample control of downstream constraint violations]
\label{prop:decision_crc}
Assume that, conditional on the training and validation data, the calibration observations and a new test observation are exchangeable,
that the trained score model $s_{\hat\vtheta}$ and threshold grid $\mathcal{T}$ are fixed independently of the calibration observations,
and that the endpoint condition \eqref{eq:crc_endpoint} holds.
For $\varepsilon\in[1/(n^{\textrm{cal}}+1),1)$, let $\hat\rho$ be the smallest grid threshold satisfying \eqref{eq:grid_crc_condition} for \eqref{eq:corrected_violation_loss}.
Then,
\begin{equation}
\label{eq:decision_crc_guarantee}
\mathbb{E}
\left[
\ell_{n^{\textrm{cal}}+1}^{\textrm{viol}}(\hat\rho)
\right]
\leq
\varepsilon,
\end{equation}
where the expectation is over calibration and test observations, conditional on the training and validation data.
\end{proposition}

The proof of
Proposition~\ref{prop:decision_crc} is given in the Appendix~\ref{app:proof_decision_crc}.

Algorithm~\ref{alg:grid_rho} summarizes the full calibration procedure.
The inputs are the calibration dataset $\mathcal D^{\textrm{cal}}$, the trained score model $s_{\hat\vtheta}$,
the ordered threshold grid $\mathcal T$, and the target violation tolerance $\varepsilon$.
First, we loop over $\mathcal{T}$ and calibration observations, construct the corresponding uncertainty set, solve the robust DCOPF problem \eqref{dcopf},
and evaluate the raw constraint violation loss (lines~\ref{alg:grid_eval_start}--\ref{alg:grid_eval_end}).
Next, for each calibration observation, we move backwards over the ordered grid $\mathcal{T}$ and recursively compute the monotone upper envelope in \eqref{eq:corrected_violation_loss}
(lines~\ref{alg:grid_corr_start}--\ref{alg:grid_corr_end}).
We then compute the corrected empirical risk $\hat R^\uparrow(\rho^j)$ at every candidate threshold
(line~\ref{alg:grid_risk_start}).
Finally, we identify all thresholds satisfying the finite-sample CRC criterion \eqref{eq:grid_crc_condition} and return the smallest feasible candidate (lines~\ref{alg:grid_feas}--\ref{alg:grid_return}).

\begin{algorithm}[tb!]
\caption{\texttt{Decision-based calibration}}
\label{alg:grid_rho}
  \textbf{Input:} Calibration dataset
  $\mathcal{D}^{\textrm{cal}}$,
  trained score model $s_{\hat\vtheta}$,
  ordered threshold grid
  $\mathcal{T}=\{\rho^1<\cdots<\rho^M\}$,
  target loss tolerance $\varepsilon$.
  \\
  \textbf{Output:} Calibrated threshold $\hat\rho$
  \begin{algorithmic}[1]
  \FOR{$j=1,\dots,M$}
    \FOR{$i=1,\dots,n^{\textrm{cal}}$}
          \STATE \label{alg:grid_eval_start}
          Construct $\mathcal{V}(\vx_i;\rho^j)$.
          \STATE Solve \eqref{dcopf} and obtain
          $\vz^\star(\vx_i;\rho^j)$.
          \STATE Evaluate
          $\ell_i^{\textrm{viol}}(\rho^j)
          =
          \mathbb{I}\Big\{
          g^{\textrm{DC}}
          \big(
          \vz^\star(\vx_i;\rho^j),
          \vxi_i
          \big)>0
          \Big\}$.
          \label{alg:grid_eval_end}
      \ENDFOR
  \ENDFOR

  \FOR{$i=1,\dots,n^{\textrm{cal}}$}
      \STATE \label{alg:grid_corr_start}
      Set
      $\ell_i^{\textrm{viol},\uparrow}(\rho^M)
      \leftarrow
      \ell_i^{\textrm{viol}}(\rho^M)$.
      \FOR{$j=M-1,\dots,1$}
          \STATE $\ell_i^{\textrm{viol},\uparrow}(\rho^j)
          \leftarrow
          \max\left\{
          \ell_i^{\textrm{viol}}(\rho^j), \ell_i^{\textrm{viol},\uparrow}(\rho^{j+1})
          \right\}$.
      \ENDFOR \label{alg:grid_corr_end}
  \ENDFOR

  \FOR{$j=1,\dots,M$}
      \STATE \label{alg:grid_risk_start}
      Compute
      $\hat R^\uparrow(\rho^j)
      =
      \frac{1}{n^{\textrm{cal}}}
      \sum_{i\in[n^{\textrm{cal}}]}
      \ell_i^{\textrm{viol},\uparrow}(\rho^j)$.
      \label{alg:grid_risk_end}
  \ENDFOR

  \STATE \label{alg:grid_feas}
  Define the feasible candidate set
  $\mathcal{T}^{\textrm{feas}}
  \leftarrow
  \left\{
  \rho^j\in\mathcal{T}:
  \hat R^\uparrow(\rho^j)
  \leq
  \varepsilon
  -
  \frac{1-\varepsilon}{n^{\textrm{cal}}}
  \right\}$.

  \STATE \label{alg:grid_return}
  Return
  $\hat\rho
  \leftarrow
  \min\mathcal{T}^{\textrm{feas}}$.

  \end{algorithmic}
\end{algorithm}

Algorithm~\ref{alg:grid_rho} requires solving at most $M n^{\textrm{cal}}$ robust DCOPF instances per trained model,
which can be trivially parallelized across observations and thresholds.
Loss evaluations can also be reused across target tolerances without solving further problem instances.
A coarse grid (small $M$) restricts threshold resolution and may require a larger $\rho$ to satisfy the CRC criterion, enlarging the uncertainty set and potentially increasing operating costs.
Nonetheless, Proposition~\ref{prop:decision_crc} remains valid under the stated assumptions.
Increasing $M$ provides more candidate thresholds at increased computational cost, but additional evaluations can reveal violations that increase the monotone upper envelope.
Therefore, a finer grid does not necessarily reduce the selected threshold.
Grid resolution balances offline computational effort against the resolution of threshold selection.
For larger downstream ARO problems, decomposition or iterative solution methods could further reduce the cost of each calibration solve.

\section{Experimental Setup and Results}
\label{sec:results}

In this section, we describe our experimental setup (Subsection~\ref{sec:exp_setup}),
and present results for an illustrative (Subsection~\ref{sec:results_5bus}) and a realistic test system (Subsection~\ref{sec:results_rts}).

\subsection{Experimental Setup}
\label{sec:exp_setup}

\paragraph{Models and metrics.}

We evaluate the reliability--cost trade-off and separately assess the effects of set shape and size calibration.
For the scores in \eqref{eq:norm_score}, we consider $p\in\{1,2,\infty\}$ and, with a slight abuse of notation, use $p=1+\infty$ for the sum of the $1$-norm and $\infty$-norm scores in \eqref{eq:sum_scores_l1inf}.
For each geometric form, we compare two shape specifications that share the same point predictor:
\begin{itemize}
    \item \texttt{CTX}: the learned contextual transformation $\mL_{\hat\vphi}(\vx)$.
    \item \texttt{ECOV}: the same point predictor $\vf_{\hat\veta}(\vx)$ with a fixed Cholesky factor $\hat{\mL}^{\textrm{emp}}$ of the empirical covariance of validation residuals. For $p\in\{1,2,\infty\}$, the corresponding score is  $ s^{\textrm{ECOV}}_{p}(\rvx,\rvy) = \left\| (\hat{\mL}^{\textrm{emp}})^{-1}
    \big(\rvy-\vf_{\hat\veta}(\rvx)\big)
    \right\|_p$;
    for $p=1+\infty$, the score is defined analogously.
\end{itemize}
We compare three calibration procedures:
\begin{itemize}
\item \texttt{CRC}: decision-based calibration with the monotonicity correction in Algorithm~\ref{alg:grid_rho}.
\item \texttt{CP}: coverage-based calibration using \eqref{eq:miscoverage_loss}.
\item \texttt{RAW}: decision-based calibration using the uncorrected loss \eqref{eq:violation_indicator}.
This procedure does not carry the finite-sample guarantee of Proposition~\ref{prop:decision_crc}.
\end{itemize}
Configurations are denoted \texttt{shape-$p$-calibration}, e.g., \texttt{CTX-2-CRC} denotes the proposed contextual $2$-norm score with decision-based calibration, whereas \texttt{ECOV-2-CP} denotes the empirical-covariance $2$-norm score with coverage-based calibration, and so forth.
For deployment, the geometric form can be selected using held-out validation data.
We assess performance for tolerances
$\varepsilon\in\{0.30,0.25,0.20,0.15,0.10,0.05\}$.
We report constraint satisfaction, defined as the fraction of test observations for which all reserve and transmission constraints hold under the scheduled recourse policy:
$1-\frac{1}{n^{\textrm{test}}}\sum_{i=1}^{n^{\textrm{test}}}
\mathbb I\!\left\{g^{\textrm{DC}}\bigl(\vz^\star(\vx_i;\hat\rho),\vxi_i\bigr)>0\right\}.$
Each observation counts once, regardless of how many constraints it violates or whether the forecast error falls inside the prediction set.
We further report the average corresponding operating cost of \eqref{dcopf} and predictive coverage, defined as the fraction of test observations inside the set.
Full results are provided in Appendix~\ref{app:full_results}.

\paragraph{Hyperparameters and model training.}
All NN models (for $\vf_{\veta}$ and $\mL_{\vphi}$) use three layers with 50 nodes per layer and ReLU hidden-layer activations.
We train with a learning rate of $10^{-3}$, $\gamma = 0.1$, and a batch size of 512 using a mini-batch gradient-based optimizer.
During the training of $\mL_{\vphi}$, we observed loss spikes and applied a small gradient clipping that resolved the issue.
We train using early stopping with a validation set ($15\%$ of training observations).
Following \cite{braun2025minimum}, we \emph{(i)} first pre-train $\vf_\veta$ by minimizing the MSE, \emph{(ii)} freeze $\vf_\veta$ and train $\mL_\vphi$ with the NLL, and \emph{(iii)} continue gradient-based training of both networks jointly using \eqref{eq:final_train_obj}.

\paragraph{End-to-end benchmarks.}
For the 5-bus system, we further consider CVaR-based end-to-end uncertainty-set learning, adapting the formulations of \citet{mieth2023prescribed,wang2023learning}.
We consider two end-to-end benchmarks: a diagonal box obtained from the $\infty$-norm score (\texttt{E2E-BOX}) and a $2$-norm score with the full contextual transformation (\texttt{E2E-2}).
We initialize each benchmark with the likelihood-trained score model and keep the point predictor fixed.
For each reliability target, we train the shape transformation through the robust DCOPF using operating cost and a CVaR-based penalty for constraint violations.
We use the training settings above and increase the penalty weight from an initial value of $500$ during training.

\paragraph{Operational setting.}
We examine two systems: the modified 5-bus system \citep{mieth2023prescribed} with 2 wind power plants
and the 73-bus, RTS-GMLC 2019 System (henceforth, we refer to it simply as RTS) \citep{barrows2019ieee} that includes 120 transmission lines, 73 conventional generators, and 4 wind power plants (we omit solar plants and assume hydro generation is fixed at the day-ahead schedule).
For the RTS system, we allow only fast gas turbines to provide reserves.
We use wind power data from the RTS system for both test systems, retaining and appropriately scaling the first two wind power series for the 5-bus system.
Consistent with the operational process outlined in Fig.~\ref{fig:operating_timeline}, each reserve schedule is computed 15 minutes ahead using the information available at the scheduling stage.
For the 5-bus system, demand is fixed, whereas for the RTS system, demand varies across scheduling instances according to the load profile but is treated as deterministic, i.e., known at the scheduling stage.
Wind power represents the uncertainty in both cases.
The feature vector $\rvx$ includes the three most recent production measurements,
following the autocorrelation analysis of \citet{stratigakos2025learning}.
We also consider an additional feature: an hourly wind power forecast issued several hours before operations, using the latest weather forecast available to SOs \citep{carmona2024_joint_granular_model}.
For both systems, we train new multivariate forecasting models.
The data span the first five months of 2020.
We use approximately $5{,}000$ observations for training and validation, $1{,}500$ for calibration, and $4{,}400$ for testing and evaluation.
The last $15\%$ of the training data is reserved for validation and hyperparameter tuning.
We compare point forecast performance with that reported by \citet{stratigakos2025learning} and obtain comparable accuracy.

\paragraph{Calibration details.}
Preliminary analysis indicated that, for a small number of RTS instances with high wind production, \eqref{dcopf} is infeasible even when $\rho=0$.
To address this, for these instances, we replace $\hat\vy$ in the nominal balance \eqref{dcopf_c1} and flow constraints \eqref{dcopf_c3}-\eqref{dcopf_c4} with a scheduled renewable injection $\vzero \leq \vq \leq \hat\vy$.
We apply the scheduled curtailment amount  $(\hat\vy-\vq)$ to realized production, hence, the realized injection becomes $\vy-(\hat\vy-\vq)=\vq+(\vy-\hat\vy)=\vq+\vxi.$
Thus, forecast errors remain defined relative to $\hat\vy$, and the recourse constraints retain their form.
We further verified that realized production exceeds scheduled curtailment in these instances.
{We note that reported operating costs include the energy-dispatch consequences of scheduled curtailment but exclude explicit curtailment penalties.}
Moreover, for both systems, \eqref{dcopf} may also become infeasible for large thresholds.
Therefore, we further add non-negative slacks to the right-hand side of the robust constraints \eqref{dcopf_c8}-\eqref{dcopf_c11}.
We heavily penalize feasibility slacks in the objective and exclude their penalties from reported operating costs.
We note that positive feasibility slacks relax the robust-feasibility assumption used in Section~\ref{sec:calibration} and, consequently, may invalidate Proposition~\ref{prop:decision_crc};
however, their activation is limited to only a small number of RTS instances and occurs only for larger thresholds, which also tend to yield overly conservative decisions.
In both systems, violations are evaluated against the original constraints without slacks at the realized forecast error;
hence, slack activation does not necessarily imply a realized violation.
For both systems, we use a grid of $M=40$ thresholds with $\rho^1=0$.
The upper endpoint $\rho^M$ is selected from preliminary coverage-calibration runs on training and validation data, and the grid is then fixed before calibration.
The same grid evaluations are reused across all target tolerances $\varepsilon$.
The finite-sample CRC guarantee assumes exchangeability of the calibration and test observations conditional on the training and validation data, as stated in Proposition~\ref{prop:decision_crc}.
We assess out-of-sample constraint satisfaction empirically on the chronological test split and recommend periodic recalibration in deployment.
Methods for non-exchangeable data based on block/bootstrap or reweighting schemes can be used when temporal dependence must be modeled explicitly \citep{farinhas2023non}.

\subsection{5-bus System}
\label{sec:results_5bus}

\begin{figure}[tb!]
\centerline{\includegraphics[width = \columnwidth]{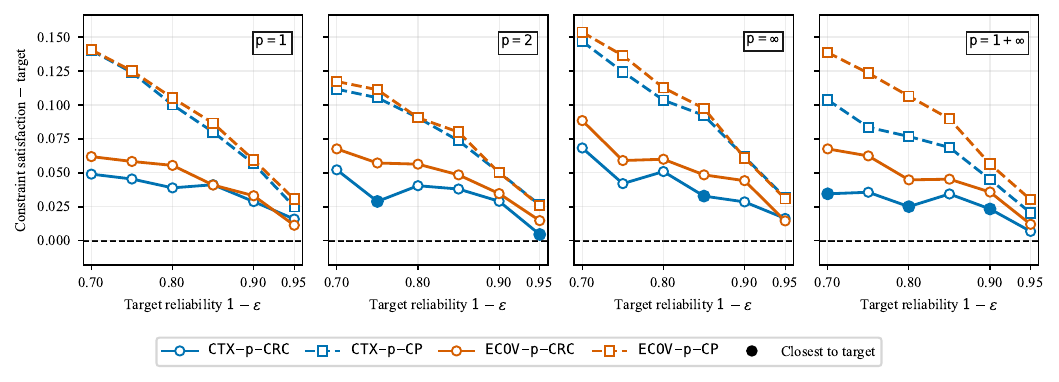}}
\caption{Difference between out-of-sample constraint satisfaction and the target reliability $1-\varepsilon$ for the 5-bus system.
The horizontal dashed line represents exact attainment of the prescribed target.
Filled markers identify the configuration closest to the target at each reliability level.}
\label{fig:violation_coverage_reliability_5bus}
\end{figure}

\paragraph{Reliability and shape--size effects.}
Figure~\ref{fig:violation_coverage_reliability_5bus} plots the difference between the achieved out-of-sample constraint satisfaction and the prescribed target reliability, with the gap above the target indicating the degree of conservativeness.
Overall, all configurations meet the prescribed reliability targets, and the gap above the target narrows as the target increases.
The ranking of configurations by proximity to the target (i.e., smaller gap is better) is largely stable across the four subplots.
Within each panel, both configurations using decision-based calibration with the monotonicity correction (\texttt{CRC}) track the targets more closely than either configuration using coverage-based calibration (\texttt{CP}).
The proposed \texttt{CTX-$p$-CRC} configurations are usually closest, with \texttt{ECOV-$p$-CRC} typically ranking second.
Moreover, the calibration procedure has a larger overall effect on conservativeness than the choice of shape transformation.
Averaged across scores, shape specifications, and reliability targets, the gap above the target falls from $8.5$ to $4.0$ percentage points when switching from \texttt{CP} to \texttt{CRC}.
Within \texttt{CRC}, the contextual transformation (\texttt{CTX}) reduces the gap by $1.3$ percentage points on average relative to the empirical covariance transformation (\texttt{ECOV}).
Among the proposed configurations, \texttt{CTX-$(1+\infty)$-CRC} tracks the targets most closely on average, followed by \texttt{CTX-2-CRC}.

\begin{figure}[tb!]
\centerline{\includegraphics[width = 3.5in]{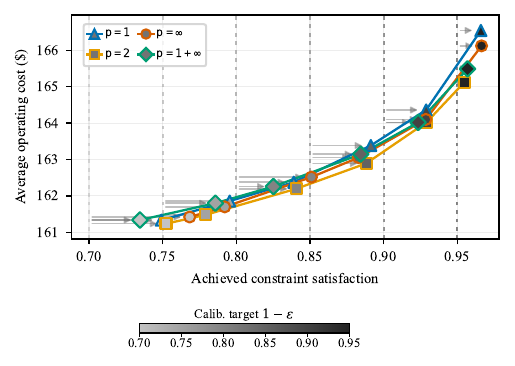}}
\caption{Reliability--cost trade-off for the proposed \texttt{CTX-$p$-CRC} configurations on the 5-bus system.
Vertical dashed lines indicate the prescribed reliability levels. Arrows show the difference between prescribed and achieved reliability.}
\label{fig:da_cost_5bus}
\end{figure}

\paragraph{Reliability--cost trade-off.}
Figure~\ref{fig:da_cost_5bus} plots average operating cost against achieved out-of-sample constraint satisfaction for the proposed \texttt{CTX-$p$-CRC} configurations.
Average operating costs rise with the reliability target for all four configurations,
with \texttt{CTX-2-CRC} and \texttt{CTX-$(1+\infty)$-CRC} offering the most favorable trade-offs.
At the $95\%$ target, for example, their markers lie just to the right of the vertical target line, indicating constraint satisfaction close to the prescribed level.
Both markers also lie below those for \texttt{CTX-1-CRC} and \texttt{CTX-$\infty$-CRC}, indicating lower operating costs.
\texttt{CTX-2-CRC} has the lowest operating cost at five of the six targets,
whereas \texttt{CTX-$(1+\infty)$-CRC} tracks the targets more closely on average at a similar cost and procures less total reserve at every target.
On average, the energy-dispatch savings under \texttt{CTX-2-CRC} more than offset its higher reserve costs.
Across the same scores, shape specifications, and reliability targets, average operating costs are approximately $0.7\%$ higher under \texttt{CP} than under \texttt{CRC}.

\paragraph{End-to-end benchmarks.}
Both end-to-end benchmarks exceed every reliability target by a wide margin.
At the $70\%$ target, for example, both achieve at least $94\%$ constraint satisfaction, compared with approximately $75\%$ for \texttt{CTX-2-CRC}.
The less costly end-to-end benchmark, \texttt{E2E-2}, has an average operating cost approximately $3.6\%$ higher than that of \texttt{CTX-2-CRC}.
Appendix~\ref{app:full_results} reports the results across all targets.

\paragraph{Effect of the monotonicity correction.}
Among the contextual configurations, the monotonicity correction changes the selected threshold only for the $2$-norm score at the $70\%$ reliability target.
In this case, the correction raises the threshold by $7.7\%$ and increases out-of-sample constraint satisfaction by about $3$ percentage points at little additional operating cost.
Both \texttt{CRC} and its uncorrected counterpart \texttt{RAW} meet all six reliability targets.
Appendix~\ref{app:full_results} reports the calibrated thresholds and test results for both procedures.

\subsection{RTS System}
\label{sec:results_rts}

\begin{figure}[tb!]
\centerline{\includegraphics[width=\columnwidth]{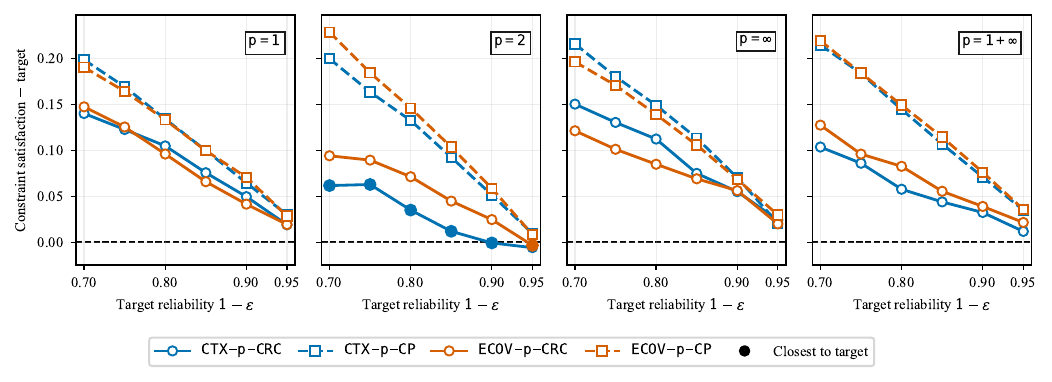}}
\caption{Difference between out-of-sample constraint satisfaction and the target reliability $1-\varepsilon$ for the RTS system.
The horizontal dashed line represents exact attainment of the prescribed target.
Filled markers identify the configuration closest to the target at each reliability level.}
\label{fig:constraint_violation_rts}
\end{figure}

\paragraph{Reliability and shape--size effects.}
As in the 5-bus system, the \texttt{CP} configurations in Figure~\ref{fig:constraint_violation_rts} meet every reliability target, often by a wide margin,
whereas the \texttt{CRC} configurations track the targets more closely in nearly all comparisons.
Namely, the mean absolute gap falls from $12.0$ percentage points under \texttt{CP} to $6.8$ under \texttt{CRC}, when averaged across scores, shape specifications, and reliability targets.
Within \texttt{CRC}, the contextual transformation reduces the mean absolute gap by $0.5$ percentage points relative to the empirical covariance transformation.
Moreover, as in the 5-bus system, the calibration procedure has a larger effect on tracking the target reliability level than the choice of shape transformation.
Among the proposed configurations, \texttt{CTX-2-CRC} tracks the targets most closely on average, followed by \texttt{CTX-$(1+\infty)$-CRC}.
\texttt{CTX-2-CRC} falls short of the $90\%$ and $95\%$ targets by less than one percentage point in each case,
whereas the other proposed configurations meet all six targets.

\begin{figure}[tb!]
\centerline{\includegraphics[width=3.5in]{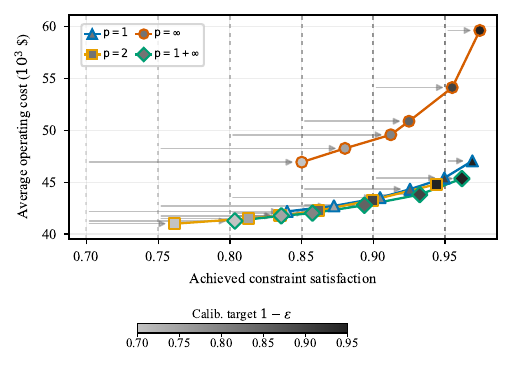}}
\caption{Reliability--cost trade-off for the proposed \texttt{CTX-$p$-CRC} configurations on the RTS system.
Vertical dashed lines indicate the prescribed reliability levels.
Arrows show the difference between prescribed and achieved reliability.}
\label{fig:da_cost_rts}
\end{figure}

\paragraph{Reliability--cost trade-off.}
Figure~\ref{fig:da_cost_rts} shows that average operating costs rise with the reliability target,
with the overall cost differences across scores being larger than in the 5-bus system.
For example, at the $95\%$ target, the marker for \texttt{CTX-2-CRC} lies just to the left of the target line, and the marker for \texttt{CTX-$(1+\infty)$-CRC} lies just to its right, and both configurations have lower operating costs than \texttt{CTX-1-CRC} and \texttt{CTX-$\infty$-CRC} at that target.
Among the proposed configurations, \texttt{CTX-2-CRC} has the lowest operating cost at every target but falls short of the two highest reliability targets,
whereas \texttt{CTX-$(1+\infty)$-CRC} meets all targets with an average cost only $0.9\%$ higher than that of \texttt{CTX-2-CRC}.
The rankings by total reserve procurement and operating cost are the same across the four scores.
Finally, the overall average operating costs are approximately $5\%$ higher under \texttt{CP} than under \texttt{CRC}, a larger difference than in the 5-bus system.

\begin{figure}[tb!]
\centerline{\includegraphics[width=0.5\columnwidth]{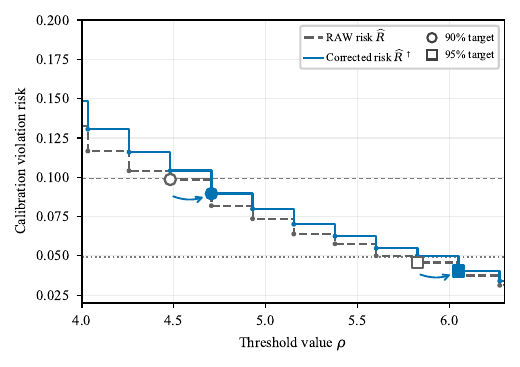}}
\caption{Raw and monotonicity-corrected empirical calibration risks for \texttt{CTX-2} in the RTS system.
Markers indicate the thresholds selected by \texttt{RAW} and \texttt{CRC} for the $90\%$ and $95\%$ reliability targets.}
\label{fig:rho_monotonicity_rts}
\end{figure}

\paragraph{Effect of the monotonicity correction.}
The monotonicity correction changes the selected threshold more often in the RTS system.
Across the four contextual scores and six reliability targets, \texttt{CRC} selects a larger threshold than \texttt{RAW} in $20$ of the $24$ cases.
Figure~\ref{fig:rho_monotonicity_rts} shows the shift for the contextual $2$-norm score at the two highest targets.
At the $95\%$ target, for example, the larger threshold increases out-of-sample constraint satisfaction by about $0.7$ percentage points, bringing it to within one percentage point of the reliability target.
Appendix~\ref{app:full_results} reports the calibrated thresholds and test results for both procedures.

\subsection{Limitations and Extensions}
\label{subsec:limitations_extensions}

The case studies use plant-level prediction sets for single-period
scheduling and evaluate one system-wide operational reliability target
at a time.
These choices motivate the following extensions.

Power-system uncertainty often has a hierarchical structure, with renewable production aggregated from individual plants to zones, balancing areas, and the full system.
Given a fixed aggregation matrix $\mH$, the joint plant-level prediction set can be mapped to the corresponding level as
$\mathcal{U}_{\mH}(\rvx;\rho) =
\{\mH\rvy:\rvy\in\mathcal{U}_{\hat\vtheta}^{\textrm{feas}}(\rvx;\rho)\}$.
The threshold $\rho$ can then be calibrated against a downstream violation loss defined at that level, following Section~\ref{sec:calibration}.
Assessing reserve deliverability still requires plant-level uncertainty and network constraints because aggregation discards the spatial distribution of forecast errors.
Reconciling sets learned separately at different levels would require
a dedicated hierarchical calibration method \citep{zhou2026hierarchical}.

Operational reliability requirements may also vary across constraint
groups. For each group $k$, we can use the same robust scheduling
solutions on the calibration grid to evaluate whether any constraint
in that group is violated. We then apply the monotonicity correction
to the group-specific loss and calibrate a threshold $\hat{\rho}_k$
at tolerance $\varepsilon_k$. Under the assumptions of
Proposition~\ref{prop:decision_crc} for each group, the common threshold
$\hat{\rho}=\max_k\hat{\rho}_k$ preserves each group's expected
violation-loss bound because each corrected loss is non-increasing in
the threshold and upper bounds the corresponding violation loss.
Nonetheless, applying the largest threshold to all groups can lead to conservative scheduling decisions when the calibrated thresholds differ substantially.

Finally, the main challenge in extending the framework to multi-period scheduling is to learn joint prediction sets for wind production trajectories that capture dependence across plants and lead times.
With the current score model, the learned transformation grows quadratically in the number of plant--lead-time pairs.
Low-rank dependence models for multivariate forecasting \citep{salinas2019high} offer a more parsimonious way to parameterize this transformation.
Intertemporal constraints such as generator ramp limits would also increase the size of the robust scheduling problems solved during
calibration, motivating approaches to speed up solution times.

\section{Conclusions}
\label{sec:conclusions}

Power system operators increasingly rely on robust optimization to balance operating costs and reliability,
but their success critically depends on the chosen uncertainty sets.
In this work, we develop a two-step framework to construct conditional multivariate prediction sets for robust power system scheduling operations.
Building on norm-based score models, we combine likelihood-based training with tractable robust reformulations and calibrate set size directly against operational constraint violations.
Applying conformal risk control to a monotone upper envelope provides finite-sample control of the expected downstream constraint-violation loss under Proposition~\ref{prop:decision_crc}.

In short-term energy and reserve scheduling experiments on the 5-bus
and RTS-GMLC systems, decision-based calibration generally achieved constraint satisfaction closer to the prescribed reliability targets and reduced average operating costs relative to coverage-based calibration.
Changing the calibration criterion had a larger effect on how closely the targets were tracked than changing the shape transformation.
Among the contextual models with decision-based calibration, the $2$-norm score achieved the lowest average operating cost in both systems and fell slightly short of two reliability targets in the RTS-GMLC system,
whereas the combined $1+\infty$ score met every reliability target in both systems at a similar average cost.

Finally, although we are primarily motivated by the operation of weather-dependent power systems, our framework could also be applied to other areas, such as portfolio optimization and multistage inventory management,
by defining the calibration loss in terms of the downstream constraints of each application.

\appendix

\section{Proof of Proposition~\ref{prop:decision_crc}}
\label{app:proof_decision_crc}

By construction, $\ell_i^{\textrm{viol},\uparrow}(\rho^j)$ is bounded by one and monotone non-increasing in $\rho^j$.
Hence, applying CRC to the corrected loss gives
\begin{equation*}
 \mathbb{E}
\left[ \ell_{n^{\textrm{cal}}+1}^{\textrm{viol},\uparrow}(\hat\rho)
\right] \leq \varepsilon
\end{equation*}
under exchangeability \citep[Theorem~1]{angelopoulos2024conformal}.
Moreover, \eqref{eq:viol_majorization} holds pointwise for every candidate threshold, so
\begin{equation*}
\ell_{n^{\textrm{cal}}+1}^{\textrm{viol}}(\hat\rho)
\leq
\ell_{n^{\textrm{cal}}+1}^{\textrm{viol},\uparrow}(\hat\rho).
\end{equation*}
Taking expectations yields \eqref{eq:decision_crc_guarantee}, which concludes the proof.

\section{Detailed Numerical Results}
\label{app:full_results}
\blue{
Tables~\ref{tab:5bus_thresholds}--\ref{tab:5bus_costs}
and Tables~\ref{tab:rts_thresholds}--\ref{tab:rts_costs}
report calibrated thresholds, constraint satisfaction,
predictive coverage, and average operating costs for
the 5-bus and RTS systems, respectively.
Column headings specify the target operational reliability
$1-\varepsilon$.
Each configuration is identified jointly by its $p$ value
and method label; for example, $p=2$ with
\texttt{CTX-CRC} corresponds to \texttt{CTX-2-CRC}.
The E2E benchmarks are included in the 5-bus performance
tables; their thresholds are omitted because
$\hat\rho=1$ is fixed during training.
Bold entries identify the smallest absolute deviation
from the column target in the constraint-satisfaction
and predictive-coverage tables, and the lowest operating
cost in the cost tables, using unrounded values.
A $\dagger$ marks configurations whose out-of-sample
constraint satisfaction falls below the operational target.
Absolute threshold values are not comparable across
score specifications because the scores operate on
different scales.
}

% 5-bus system

\begin{table*}[tb!]
\centering
\caption{\blue{Calibrated thresholds $\hat\rho$ for the 5-bus system.}}
\label{tab:5bus_thresholds}
\small
\setlength{\tabcolsep}{5.5pt}
\begin{tabular}{clrrrrrr}
\toprule
& & \multicolumn{6}{c}{Target reliability $1-\varepsilon$} \\
\cmidrule(lr){3-8}
$p$ & Method & 70\% & 75\% & 80\% & 85\% & 90\% & 95\% \\
\midrule
$1$ & \texttt{CTX-CRC} & 1.862 & 2.149 & 2.435 & 3.008 & 3.581 & 4.871 \\
 & \texttt{CTX-RAW} & 1.862 & 2.149 & 2.435 & 3.008 & 3.581 & 4.871 \\
 & \texttt{CTX-CP} & 2.450 & 2.773 & 3.141 & 3.622 & 4.482 & 5.587 \\
 & \texttt{ECOV-CRC} & 1.104 & 1.288 & 1.564 & 1.840 & 2.300 & 2.852 \\
 & \texttt{ECOV-CP} & 1.469 & 1.693 & 1.970 & 2.351 & 2.819 & 3.589 \\
\cmidrule(lr){1-8}
$2$ & \texttt{CTX-CRC} & 1.909 & 2.045 & 2.455 & 2.864 & 3.545 & 4.227 \\
 & \texttt{CTX-RAW} & 1.773 & 2.045 & 2.455 & 2.864 & 3.545 & 4.227 \\
 & \texttt{CTX-CP} & 2.265 & 2.564 & 2.901 & 3.367 & 4.089 & 5.318 \\
 & \texttt{ECOV-CRC} & 0.973 & 1.122 & 1.347 & 1.646 & 2.020 & 2.544 \\
 & \texttt{ECOV-CP} & 1.174 & 1.380 & 1.581 & 1.937 & 2.273 & 2.918 \\
\cmidrule(lr){1-8}
$\infty$ & \texttt{CTX-CRC} & 2.014 & 2.168 & 2.633 & 2.943 & 3.562 & 4.801 \\
 & \texttt{CTX-RAW} & 2.014 & 2.168 & 2.633 & 2.943 & 3.562 & 4.801 \\
 & \texttt{CTX-CP} & 2.589 & 2.858 & 3.218 & 3.916 & 4.651 & 6.041 \\
 & \texttt{ECOV-CRC} & 0.840 & 0.910 & 1.119 & 1.329 & 1.749 & 2.169 \\
 & \texttt{ECOV-CP} & 1.090 & 1.259 & 1.447 & 1.795 & 2.121 & 2.729 \\
\cmidrule(lr){1-8}
$1+\infty$ & \texttt{CTX-CRC} & 1.773 & 2.010 & 2.246 & 2.719 & 3.192 & 4.019 \\
 & \texttt{CTX-RAW} & 1.773 & 2.010 & 2.246 & 2.719 & 3.192 & 4.019 \\
 & \texttt{CTX-CP} & 2.110 & 2.315 & 2.656 & 3.117 & 3.697 & 4.610 \\
 & \texttt{ECOV-CRC} & 1.953 & 2.278 & 2.603 & 3.254 & 4.068 & 5.044 \\
 & \texttt{ECOV-CP} & 2.557 & 2.944 & 3.431 & 4.180 & 4.863 & 6.346 \\
\bottomrule
\end{tabular}
\vspace{1mm}

\parbox{0.96\textwidth}{\footnotesize \emph{Note:} Absolute values are not comparable across score functions because the scores operate on different scales.}
\end{table*}

\begin{table*}[tb!]
\centering
\caption{\blue{Out-of-sample constraint satisfaction (\%) for the 5-bus system.}}
\label{tab:5bus_constraint_satisfaction}
\small
\setlength{\tabcolsep}{5.5pt}
\begin{tabular}{clrrrrrr}
\toprule
& & \multicolumn{6}{c}{Target reliability $1-\varepsilon$} \\
\cmidrule(lr){3-8}
$p$ & Method & 70\% & 75\% & 80\% & 85\% & 90\% & 95\% \\
\midrule
$1$ & \texttt{CTX-CRC} & 74.91 & 79.55 & 83.90 & 89.13 & 92.89 & 96.60 \\
 & \texttt{CTX-RAW} & 74.91 & 79.55 & 83.90 & 89.13 & 92.89 & 96.60 \\
 & \texttt{CTX-CP} & 84.01 & 87.39 & 90.01 & 93.00 & 95.67 & 97.51 \\
 & \texttt{ECOV-CRC} & 76.20 & 80.84 & 85.55 & 89.11 & 93.32 & 96.15 \\
 & \texttt{ECOV-CP} & 84.06 & 87.50 & 90.51 & 93.66 & 95.95 & 98.08 \\
\cmidrule(lr){1-8}
$2$ & \texttt{CTX-CRC} & 75.23 & \textbf{77.90} & 84.06 & 88.81 & 92.91 & \textbf{95.47} \\
 & \texttt{CTX-RAW} & \textbf{72.31} & \textbf{77.90} & 84.06 & 88.81 & 92.91 & \textbf{95.47} \\
 & \texttt{CTX-CP} & 81.16 & 85.53 & 89.06 & 92.37 & 95.04 & 97.64 \\
 & \texttt{ECOV-CRC} & 76.77 & 80.73 & 85.64 & 89.86 & 93.48 & 96.49 \\
 & \texttt{ECOV-CP} & 81.73 & 86.12 & 89.06 & 93.00 & 95.02 & 97.60 \\
\cmidrule(lr){1-8}
$\infty$ & \texttt{CTX-CRC} & 76.83 & 79.21 & 85.10 & \textbf{88.29} & 92.87 & 96.65 \\
 & \texttt{CTX-RAW} & 76.83 & 79.21 & 85.10 & \textbf{88.29} & 92.87 & 96.65 \\
 & \texttt{CTX-CP} & 84.67 & 87.43 & 90.35 & 94.23 & 96.17 & 98.19 \\
 & \texttt{ECOV-CRC} & 78.85 & 80.91 & 86.01 & 89.86 & 94.43 & 96.47 \\
 & \texttt{ECOV-CP} & 85.35 & 88.63 & 91.26 & 94.75 & 96.11 & 98.10 \\
\cmidrule(lr){1-8}
$1+\infty$ & \texttt{CTX-CRC} & 73.46 & 78.58 & \textbf{82.52} & 88.45 & \textbf{92.35} & 95.70 \\
 & \texttt{CTX-RAW} & 73.46 & 78.58 & \textbf{82.52} & 88.45 & \textbf{92.35} & 95.70 \\
 & \texttt{CTX-CP} & 80.37 & 83.36 & 87.68 & 91.87 & 94.52 & 97.06 \\
 & \texttt{ECOV-CRC} & 76.77 & 81.25 & 84.49 & 89.54 & 93.59 & 96.22 \\
 & \texttt{ECOV-CP} & 83.85 & 87.34 & 90.65 & 93.98 & 95.65 & 98.03 \\
\cmidrule(lr){1-8}
\multicolumn{1}{c}{--} & \texttt{E2E-BOX} & 98.60 & 98.80 & 99.12 & 99.34 & 99.48 & 99.66 \\
\multicolumn{1}{c}{--} & \texttt{E2E-2} & 94.38 & 96.88 & 97.98 & 98.64 & 99.21 & 99.57 \\
\bottomrule
\end{tabular}
\vspace{1mm}

\parbox{0.96\textwidth}{\footnotesize \emph{Note:} Bold entries are closest to the target within each column.}
\end{table*}

\begin{table*}[tb!]
\centering
\caption{\blue{Out-of-sample predictive coverage (\%) for the 5-bus system.}}
\label{tab:5bus_predictive_coverage}
\small
\setlength{\tabcolsep}{5.5pt}
\begin{tabular}{clrrrrrr}
\toprule
& & \multicolumn{6}{c}{Target reliability $1-\varepsilon$} \\
\cmidrule(lr){3-8}
$p$ & Method & 70\% & 75\% & 80\% & 85\% & 90\% & 95\% \\
\midrule
$1$ & \texttt{CTX-CRC} & 59.94 & 66.53 & 72.76 & 81.34 & 86.78 & 93.95 \\
 & \texttt{CTX-RAW} & 59.94 & 66.53 & 72.76 & 81.34 & 86.78 & 93.95 \\
 & \texttt{CTX-CP} & 73.10 & 78.01 & 82.81 & 87.25 & 92.14 & \textbf{95.52} \\
 & \texttt{ECOV-CRC} & 65.90 & 71.76 & 77.94 & 82.70 & 88.20 & 92.53 \\
 & \texttt{ECOV-CP} & 76.20 & 80.37 & 84.62 & 89.04 & 92.35 & 95.99 \\
\cmidrule(lr){1-8}
$2$ & \texttt{CTX-CRC} & 64.42 & 67.55 & 75.41 & 81.41 & 88.25 & 92.44 \\
 & \texttt{CTX-RAW} & 61.07 & 67.55 & 75.41 & 81.41 & 88.25 & 92.44 \\
 & \texttt{CTX-CP} & 71.72 & 77.11 & 82.11 & 86.93 & 91.62 & 95.61 \\
 & \texttt{ECOV-CRC} & \textbf{69.29} & \textbf{74.28} & \textbf{80.16} & \textbf{85.71} & 90.58 & 94.20 \\
 & \texttt{ECOV-CP} & 75.72 & 80.75 & 84.60 & 89.72 & 92.69 & 95.86 \\
\cmidrule(lr){1-8}
$\infty$ & \texttt{CTX-CRC} & 64.58 & 67.89 & 75.88 & 80.03 & 86.59 & 93.84 \\
 & \texttt{CTX-RAW} & 64.58 & 67.89 & 75.88 & 80.03 & 86.59 & 93.84 \\
 & \texttt{CTX-CP} & 75.09 & 79.14 & 83.11 & 89.47 & 93.12 & 96.65 \\
 & \texttt{ECOV-CRC} & 66.85 & 69.90 & 76.95 & 82.27 & \textbf{89.65} & 93.07 \\
 & \texttt{ECOV-CP} & 76.06 & 80.75 & 84.74 & 90.04 & 92.75 & 95.90 \\
\cmidrule(lr){1-8}
$1+\infty$ & \texttt{CTX-CRC} & 64.74 & 70.95 & 75.79 & 83.29 & 88.18 & 93.68 \\
 & \texttt{CTX-RAW} & 64.74 & 70.95 & 75.79 & 83.29 & 88.18 & 93.68 \\
 & \texttt{CTX-CP} & 72.94 & 76.97 & 82.00 & 87.64 & 92.03 & 95.67 \\
 & \texttt{ECOV-CRC} & 66.76 & 72.44 & 76.45 & 83.27 & 88.95 & 93.05 \\
 & \texttt{ECOV-CP} & 75.82 & 80.46 & 84.51 & 89.67 & 92.39 & 95.90 \\
\cmidrule(lr){1-8}
\multicolumn{1}{c}{--} & \texttt{E2E-BOX} & 96.69 & 97.26 & 97.80 & 98.39 & 98.91 & 99.34 \\
\multicolumn{1}{c}{--} & \texttt{E2E-2} & 79.87 & 87.73 & 93.07 & 95.70 & 97.62 & 98.75 \\
\bottomrule
\end{tabular}
\vspace{1mm}

\parbox{0.96\textwidth}{\footnotesize \emph{Note:} Bold entries are closest to the target within each column.}
\end{table*}

\begin{table*}[tb!]
\centering
\caption{\blue{Average operating cost (\$) per 15-minute period for the 5-bus system.}}
\label{tab:5bus_costs}
\small
\setlength{\tabcolsep}{5.5pt}
\begin{tabular}{clrrrrrr}
\toprule
& & \multicolumn{6}{c}{Target reliability $1-\varepsilon$} \\
\cmidrule(lr){3-8}
$p$ & Method & 70\% & 75\% & 80\% & 85\% & 90\% & 95\% \\
\midrule
$1$ & \texttt{CTX-CRC} & 161.33 & 161.86 & 162.38 & 163.38 & 164.36 & 166.55 \\
 & \texttt{CTX-RAW} & 161.33 & 161.86 & 162.38 & 163.38 & 164.36 & 166.55 \\
 & \texttt{CTX-CP} & 162.40 & 162.97 & 163.61 & 164.43 & 165.87 & 167.81 \\
 & \texttt{ECOV-CRC} & 161.24 & 161.77 & 162.54 & 163.28 & 164.47 & 165.84 \\
 & \texttt{ECOV-CP} & 162.27 & 162.88 & 163.62 & 164.60 & 165.75 & 167.56 \\
\cmidrule(lr){1-8}
$2$ & \texttt{CTX-CRC} & 161.24 & \textbf{161.48} & \textbf{162.19} & \textbf{162.89} & 164.02 & \textbf{165.12} \\
 & \texttt{CTX-RAW} & \textbf{160.99} & \textbf{161.48} & \textbf{162.19} & \textbf{162.89} & 164.02 & \textbf{165.12} \\
 & \texttt{CTX-CP} & 161.86 & 162.38 & 162.95 & 163.73 & 164.90 & 166.83 \\
 & \texttt{ECOV-CRC} & 161.35 & 161.86 & 162.59 & 163.54 & 164.68 & 166.22 \\
 & \texttt{ECOV-CP} & 162.03 & 162.70 & 163.34 & 164.43 & 165.43 & 167.27 \\
\cmidrule(lr){1-8}
$\infty$ & \texttt{CTX-CRC} & 161.42 & 161.70 & 162.52 & 163.05 & 164.10 & 166.12 \\
 & \texttt{CTX-RAW} & 161.42 & 161.70 & 162.52 & 163.05 & 164.10 & 166.12 \\
 & \texttt{CTX-CP} & 162.44 & 162.91 & 163.52 & 164.69 & 165.88 & 168.07 \\
 & \texttt{ECOV-CRC} & 161.59 & 161.88 & 162.73 & 163.56 & 165.16 & 166.68 \\
 & \texttt{ECOV-CP} & 162.62 & 163.29 & 164.02 & 165.33 & 166.51 & 168.63 \\
\cmidrule(lr){1-8}
$1+\infty$ & \texttt{CTX-CRC} & 161.34 & 161.80 & 162.26 & 163.15 & \textbf{164.02} & 165.49 \\
 & \texttt{CTX-RAW} & 161.34 & 161.80 & 162.26 & 163.15 & \textbf{164.02} & 165.49 \\
 & \texttt{CTX-CP} & 161.99 & 162.39 & 163.03 & 163.88 & 164.92 & 166.54 \\
 & \texttt{ECOV-CRC} & 161.23 & 161.75 & 162.26 & 163.26 & 164.46 & 165.83 \\
 & \texttt{ECOV-CP} & 162.19 & 162.79 & 163.53 & 164.62 & 165.58 & 167.58 \\
\cmidrule(lr){1-8}
\multicolumn{1}{c}{--} & \texttt{E2E-BOX} & 168.14 & 168.77 & 169.54 & 170.56 & 172.15 & 174.85 \\
\multicolumn{1}{c}{--} & \texttt{E2E-2} & 166.58 & 166.83 & 167.41 & 168.38 & 169.73 & 173.14 \\
\bottomrule
\end{tabular}
\vspace{1mm}

\parbox{0.96\textwidth}{\footnotesize \emph{Note:} Bold entries are the lowest cost within each column.}
\end{table*}

% RTS system

\begin{table*}[tb!]
\centering
\caption{\blue{Calibrated thresholds $\hat\rho$ for the RTS system.}}
\label{tab:rts_thresholds}
\small
\setlength{\tabcolsep}{5.5pt}
\begin{tabular}{clrrrrrr}
\toprule
& & \multicolumn{6}{c}{Target reliability $1-\varepsilon$} \\
\cmidrule(lr){3-8}
$p$ & Method & 70\% & 75\% & 80\% & 85\% & 90\% & 95\% \\
\midrule
$1$ & \texttt{CTX-CRC} & 3.774 & 4.246 & 4.954 & 5.661 & 6.605 & 8.020 \\
 & \texttt{CTX-RAW} & 3.538 & 4.010 & 4.718 & 5.425 & 6.369 & 7.784 \\
 & \texttt{CTX-CP} & 4.791 & 5.385 & 5.929 & 6.636 & 7.614 & 9.200 \\
 & \texttt{ECOV-CRC} & 2.646 & 2.977 & 3.473 & 3.969 & 4.465 & 5.457 \\
 & \texttt{ECOV-CP} & 3.495 & 3.852 & 4.280 & 4.819 & 5.475 & 6.450 \\
\cmidrule(lr){1-8}
$2$ & \texttt{CTX-CRC} & 2.688 & 3.136 & 3.360 & 3.808 & 4.705\textsuperscript{\(\dagger\)} & 6.049\textsuperscript{\(\dagger\)} \\
 & \texttt{CTX-RAW} & 2.688 & 2.912 & 3.360 & 3.808 & 4.481\textsuperscript{\(\dagger\)} & 5.825\textsuperscript{\(\dagger\)} \\
 & \texttt{CTX-CP} & 4.523 & 4.947 & 5.578 & 6.156 & 7.181 & 8.737 \\
 & \texttt{ECOV-CRC} & 1.217 & 1.420 & 1.623 & 1.826 & 2.130 & 2.637\textsuperscript{\(\dagger\)} \\
 & \texttt{ECOV-CP} & 2.146 & 2.345 & 2.598 & 2.872 & 3.328 & 3.956 \\
\cmidrule(lr){1-8}
$\infty$ & \texttt{CTX-CRC} & 2.770 & 3.196 & 3.622 & 4.048 & 5.113 & 6.818 \\
 & \texttt{CTX-RAW} & 2.130 & 2.557 & 2.983 & 3.409 & 4.261 & 5.965 \\
 & \texttt{CTX-CP} & 3.830 & 4.205 & 4.715 & 5.479 & 6.432 & 8.309 \\
 & \texttt{ECOV-CRC} & 1.201 & 1.372 & 1.544 & 1.801 & 2.230 & 3.002 \\
 & \texttt{ECOV-CP} & 1.655 & 1.842 & 2.056 & 2.294 & 2.747 & 3.345 \\
\cmidrule(lr){1-8}
$1+\infty$ & \texttt{CTX-CRC} & 3.385 & 3.906 & 4.166 & 4.948 & 5.989 & 7.552 \\
 & \texttt{CTX-RAW} & 3.125 & 3.646 & 4.166 & 4.687 & 5.729 & 7.291 \\
 & \texttt{CTX-CP} & 5.501 & 6.016 & 6.676 & 7.330 & 8.300 & 10.156 \\
 & \texttt{ECOV-CRC} & 3.454 & 3.701 & 4.194 & 4.934 & 5.921 & 7.648 \\
 & \texttt{ECOV-CP} & 5.284 & 5.701 & 6.353 & 7.142 & 8.177 & 9.622 \\
\bottomrule
\end{tabular}
\vspace{1mm}

\parbox{0.96\textwidth}{\footnotesize \emph{Note:} Absolute values are not comparable across score functions because the scores operate on different scales.}
\end{table*}

\begin{table*}[tb!]
\centering
\caption{\blue{Out-of-sample constraint satisfaction (\%) for the RTS system.}}
\label{tab:rts_constraint_satisfaction}
\small
\setlength{\tabcolsep}{5.5pt}
\begin{tabular}{clrrrrrr}
\toprule
& & \multicolumn{6}{c}{Target reliability $1-\varepsilon$} \\
\cmidrule(lr){3-8}
$p$ & Method & 70\% & 75\% & 80\% & 85\% & 90\% & 95\% \\
\midrule
$1$ & \texttt{CTX-CRC} & 83.99 & 87.25 & 90.47 & 92.55 & 94.95 & 96.90 \\
 & \texttt{CTX-RAW} & 82.88 & 85.69 & 89.29 & 92.03 & 94.57 & 96.69 \\
 & \texttt{CTX-CP} & 89.83 & 91.89 & 93.39 & 94.97 & 96.44 & 97.98 \\
 & \texttt{ECOV-CRC} & 84.71 & 87.52 & 89.58 & 91.58 & 94.13 & 96.92 \\
 & \texttt{ECOV-CP} & 89.02 & 91.39 & 93.30 & 94.95 & 97.06 & 97.85 \\
\cmidrule(lr){1-8}
$2$ & \texttt{CTX-CRC} & \textbf{76.15} & 81.27 & \textbf{83.49} & \textbf{86.19} & \textbf{89.92}\textsuperscript{\(\dagger\)} & 94.41\textsuperscript{\(\dagger\)} \\
 & \texttt{CTX-RAW} & \textbf{76.15} & \textbf{79.33} & \textbf{83.49} & \textbf{86.19} & 89.90\textsuperscript{\(\dagger\)} & 93.68\textsuperscript{\(\dagger\)} \\
 & \texttt{CTX-CP} & 89.95 & 91.33 & 93.23 & 94.23 & 95.13 & 95.95 \\
 & \texttt{ECOV-CRC} & 79.39 & 83.92 & 87.12 & 89.47 & 92.46 & \textbf{94.66}\textsuperscript{\(\dagger\)} \\
 & \texttt{ECOV-CP} & 92.82 & 93.41 & 94.57 & 95.34 & 95.86 & 95.83 \\
\cmidrule(lr){1-8}
$\infty$ & \texttt{CTX-CRC} & 85.01 & 88.02 & 91.21 & 92.48 & 95.49 & 97.42 \\
 & \texttt{CTX-RAW} & 79.53 & 83.81 & 86.53 & 89.79 & 93.30 & 96.74 \\
 & \texttt{CTX-CP} & 91.55 & 93.00 & 94.88 & 96.31 & 97.01 & 97.08 \\
 & \texttt{ECOV-CRC} & 82.11 & 85.10 & 88.47 & 91.89 & 95.61 & 96.99 \\
 & \texttt{ECOV-CP} & 89.61 & 92.05 & 93.89 & 95.56 & 96.83 & 97.96 \\
\cmidrule(lr){1-8}
$1+\infty$ & \texttt{CTX-CRC} & 80.34 & 83.58 & 85.76 & 89.38 & 93.23 & 96.17 \\
 & \texttt{CTX-RAW} & 77.76 & 81.95 & 85.76 & 88.04 & 92.53 & 95.65 \\
 & \texttt{CTX-CP} & 91.44 & 93.34 & 94.43 & 95.67 & 97.08 & 98.48 \\
 & \texttt{ECOV-CRC} & 82.72 & 84.58 & 88.25 & 90.53 & 93.89 & 97.15 \\
 & \texttt{ECOV-CP} & 91.89 & 93.39 & 94.90 & 96.42 & 97.60 & 98.57 \\
\bottomrule
\end{tabular}
\vspace{1mm}

\parbox{0.96\textwidth}{\footnotesize \emph{Note:} Bold entries are closest to the target within each column. A $\dagger$ marks a configuration whose out-of-sample constraint satisfaction is below the prescribed target.}
\end{table*}

\begin{table*}[tb!]
\centering
\caption{\blue{Out-of-sample predictive coverage (\%) for the RTS system.}}
\label{tab:rts_predictive_coverage}
\small
\setlength{\tabcolsep}{5.5pt}
\begin{tabular}{clrrrrrr}
\toprule
& & \multicolumn{6}{c}{Target reliability $1-\varepsilon$} \\
\cmidrule(lr){3-8}
$p$ & Method & 70\% & 75\% & 80\% & 85\% & 90\% & 95\% \\
\midrule
$1$ & \texttt{CTX-CRC} & 66.24 & 72.37 & 80.21 & \textbf{85.53} & \textbf{90.13} & 94.02 \\
 & \texttt{CTX-RAW} & 62.79 & 69.41 & 77.81 & 84.08 & 89.27 & 93.68 \\
 & \texttt{CTX-CP} & 78.62 & 83.61 & 86.89 & 90.17 & 93.41 & 95.90 \\
 & \texttt{ECOV-CRC} & \textbf{67.96} & \textbf{73.14} & \textbf{79.87} & 84.33 & 87.79 & 93.14 \\
 & \texttt{ECOV-CP} & 79.98 & 83.40 & 86.59 & 89.81 & 93.18 & 96.06 \\
\cmidrule(lr){1-8}
$2$ & \texttt{CTX-CRC} & 49.00 & 58.27 & 62.82 & 68.98 & 79.39\textsuperscript{\(\dagger\)} & 88.88\textsuperscript{\(\dagger\)} \\
 & \texttt{CTX-RAW} & 49.00 & 53.74 & 62.82 & 68.98 & 77.20\textsuperscript{\(\dagger\)} & 87.18\textsuperscript{\(\dagger\)} \\
 & \texttt{CTX-CP} & 77.67 & 81.34 & 85.62 & 89.47 & 92.46 & 95.70 \\
 & \texttt{ECOV-CRC} & 53.92 & 61.41 & 68.21 & 73.21 & 79.26 & 86.89\textsuperscript{\(\dagger\)} \\
 & \texttt{ECOV-CP} & 79.46 & 82.90 & 86.48 & 89.38 & 92.98 & 95.63 \\
\cmidrule(lr){1-8}
$\infty$ & \texttt{CTX-CRC} & 65.26 & 71.88 & 77.20 & 81.11 & 88.59 & 94.32 \\
 & \texttt{CTX-RAW} & 51.92 & 60.91 & 68.98 & 74.77 & 82.81 & 92.44 \\
 & \texttt{CTX-CP} & 79.10 & 82.50 & 86.53 & 90.38 & 93.59 & 96.60 \\
 & \texttt{ECOV-CRC} & 63.72 & 70.13 & 75.05 & 80.98 & 88.07 & 94.36 \\
 & \texttt{ECOV-CP} & 77.49 & 81.79 & 85.82 & 88.81 & 92.73 & 95.88 \\
\cmidrule(lr){1-8}
$1+\infty$ & \texttt{CTX-CRC} & 47.60 & 57.25 & 61.50 & 71.15 & 79.98 & 88.75 \\
 & \texttt{CTX-RAW} & 42.01 & 52.79 & 61.50 & 68.12 & 77.76 & 87.55 \\
 & \texttt{CTX-CP} & 76.25 & 80.05 & 84.49 & 87.73 & 91.06 & \textbf{95.24} \\
 & \texttt{ECOV-CRC} & 61.30 & 64.61 & 70.65 & 77.29 & 84.22 & 91.64 \\
 & \texttt{ECOV-CP} & 79.98 & 82.88 & 86.30 & 89.74 & 93.18 & 95.79 \\
\bottomrule
\end{tabular}
\vspace{1mm}

\parbox{0.96\textwidth}{\footnotesize \emph{Note:} Bold entries are closest to the target within each column. 
A $\dagger$ marks a configuration whose out-of-sample constraint satisfaction is below the prescribed target.}
\end{table*}

\begin{table*}[tb!]
\centering
\caption{\blue{Average operating cost (thousand dollars) per 15-minute period for the RTS system.}}
\label{tab:rts_costs}
\small
\setlength{\tabcolsep}{5.5pt}
\begin{tabular}{clrrrrrr}
\toprule
& & \multicolumn{6}{c}{Target reliability $1-\varepsilon$} \\
\cmidrule(lr){3-8}
$p$ & Method & 70\% & 75\% & 80\% & 85\% & 90\% & 95\% \\
\midrule
$1$ & \texttt{CTX-CRC} & 42.18 & 42.71 & 43.51 & 44.32 & 45.41 & 47.05 \\
 & \texttt{CTX-RAW} & 41.92 & 42.44 & 43.24 & 44.05 & 45.13 & 46.78 \\
 & \texttt{CTX-CP} & 43.32 & 44.00 & 44.62 & 45.44 & 46.58 & 48.43 \\
 & \texttt{ECOV-CRC} & 42.06 & 42.55 & 43.28 & 44.05 & 44.82 & 46.34 \\
 & \texttt{ECOV-CP} & 43.31 & 43.86 & 44.53 & 45.36 & 46.37 & 47.86 \\
\cmidrule(lr){1-8}
$2$ & \texttt{CTX-CRC} & \textbf{40.99} & 41.49 & \textbf{41.74} & \textbf{42.25} & 43.27\textsuperscript{\(\dagger\)} & 44.81\textsuperscript{\(\dagger\)} \\
 & \texttt{CTX-RAW} & \textbf{40.99} & \textbf{41.24} & \textbf{41.74} & \textbf{42.25} & \textbf{43.01}\textsuperscript{\(\dagger\)} & 44.55\textsuperscript{\(\dagger\)} \\
 & \texttt{CTX-CP} & 43.06 & 43.55 & 44.27 & 44.94 & 46.12 & 47.91 \\
 & \texttt{ECOV-CRC} & 41.05 & 41.54 & 42.03 & 42.51 & 43.23 & \textbf{44.49}\textsuperscript{\(\dagger\)} \\
 & \texttt{ECOV-CP} & 43.26 & 43.76 & 44.39 & 45.07 & 46.19 & 47.72 \\
\cmidrule(lr){1-8}
$\infty$ & \texttt{CTX-CRC} & 46.93 & 48.26 & 49.58 & 50.88 & 54.13 & 59.65 \\
 & \texttt{CTX-RAW} & 44.92 & 46.26 & 47.60 & 48.92 & 51.53 & 56.82 \\
 & \texttt{CTX-CP} & 50.22 & 51.36 & 52.91 & 55.27 & 58.36 & 64.96 \\
 & \texttt{ECOV-CRC} & 43.08 & 43.79 & 44.51 & 45.58 & 47.35 & 50.47 \\
 & \texttt{ECOV-CP} & 44.97 & 45.75 & 46.63 & 47.61 & 49.44 & 51.82 \\
\cmidrule(lr){1-8}
$1+\infty$ & \texttt{CTX-CRC} & 41.25 & 41.76 & 42.01 & 42.77 & 43.80 & 45.35 \\
 & \texttt{CTX-RAW} & 41.00 & 41.50 & 42.01 & 42.52 & 43.54 & 45.09 \\
 & \texttt{CTX-CP} & 43.32 & 43.83 & 44.48 & 45.13 & 46.09 & 47.94 \\
 & \texttt{ECOV-CRC} & 41.71 & 41.97 & 42.47 & 43.23 & 44.29 & 46.13 \\
 & \texttt{ECOV-CP} & 43.60 & 44.05 & 44.75 & 45.59 & 46.68 & 48.20 \\
\bottomrule
\end{tabular}
\vspace{1mm}

\parbox{0.96\textwidth}{\footnotesize \emph{Note:} Bold entries are the lowest cost within each column. A $\dagger$ marks a configuration whose out-of-sample constraint satisfaction is below the prescribed target.}
\end{table*}

\bibliographystyle{informs2014} % outcomment this and next line in Case 1
\bibliography{references} % if more than one, comma separated

\end{document}